\documentclass{article}
\usepackage{amsmath}
\usepackage{mathrsfs}
\usepackage{txfonts}
\usepackage{bm}
\usepackage{latexsym,amsfonts,amssymb}
\begin{document}
\setcounter{page}{1}
\newtheorem{thm}{Theorem}[section]
\newtheorem{fthm}[thm]{Fundamental Theorem}
\newtheorem{dfn}[thm]{Definition}
\newtheorem{rem}[thm]{Remark}
\newtheorem{lem}[thm]{Lemma}
\newtheorem{cor}[thm]{Corollary}
\newtheorem{exa}[thm]{Example}
\newtheorem{pro}[thm]{Proposition}
\newtheorem{prob}[thm]{Problem}

\newtheorem{theorem}{Theorem}[section]
\newtheorem{definition}[theorem]{Definition}
\newtheorem{remark}[theorem]{Remark}
\newtheorem{lemma}[theorem]{Lemma}
\newtheorem{corollary}[theorem]{Corollary}
\newtheorem{example}[theorem]{Example}
\newtheorem{proposition}[theorem]{Proposition}
\newtheorem{con}[thm]{Conjecture}
\newtheorem{ob}[thm]{Observation}

\renewcommand{\theequation}{\thesection.\arabic{equation}}
\renewcommand{\thefootnote}{\fnsymbol{footnote}}

\renewcommand{\thefootnote}{\fnsymbol{footnote}}
\newcommand{\qed}{\hfill\Box\medskip}
\newcommand{\proof}{{\it Proof.\quad}}

\newcommand{\rmnum}[1]{\romannumeral #1}
\renewcommand{\abovewithdelims}[2]{%
\genfrac{[}{]}{0pt}{}{#1}{#2}}

\renewcommand{\thefootnote}{\fnsymbol{footnote}}

\title{\bf Vector spaces and Grassmann graphs \\over residue class rings\footnote{L. Huang is supported by NSFC (11371072), B. Lv is supported by NSFC (11501036, 11301270), K. Wang is supported by NSFC
( 11671043, 11371204).  Supported by  the
Fundamental Research Funds for the Central University of China.}}

\author{Li-Ping Huang\textsuperscript{a}\footnote{E-mail address: lipingmath@163.com (L. Huang), bjlv@bnu.edu.cn (B. Lv), wangks@bnu.edu.cn (K. Wang)}
\quad Benjian Lv\textsuperscript{b}\quad Kaishun Wang\textsuperscript{b}\\
{\footnotesize  \textsuperscript{a} \em  School of Math. and Statis.,
Changsha University of Science and
Technology, Changsha, 410004, China}\\
{\footnotesize  \textsuperscript{b} \em Sch. Math. Sci. {\rm \&}
Lab. Math. Com. Sys., Beijing Normal University, Beijing, 100875,
China } }
 \date{}
 \maketitle

\begin{abstract}
Let $\mathbb{Z}_{p^s}$ be the residue class ring of integers modulo $p^s$, where $p$ is a prime number and $s$ is a positive integer.
 Using  matrix representation and the inner rank of a matrix, we study the intersection, join, dimension formula and dual subspaces on
 vector subspaces of $\mathbb{Z}^n_{p^s}$. Based on these results, we investigate the Grassmann graph $G_{p^s}(n,m)$ over  $\mathbb{Z}_{p^s}$.
 $G_{p^s}(n,m)$ is a connected vertex-transitive graph, and we determine its valency,  clique number and maximum cliques. Finally,
we characterize the automorphisms of $G_{p^s}(n,m)$.

\vspace{2mm}

\noindent{\bf Keywords:}  residue class ring, subspace, dimension formula, Grassmann graph,  maximum clique,  automorphism

\vspace{2mm}

\noindent{\bf 2010 AMS Classification}:  05C25, 15B33,  05C60, 05B25,  51A10

\end{abstract}

\section{Introduction}
 \setcounter{equation}{0}

\ \ \ \ \
Throughout, let $R$ be a commutative local ring and  $R^*$  the set of all units of $R$.  For a subset $S$ of $R$, let $S^{m\times n}$ be the set of all $m\times n$ matrices over $S$,
and $S^n=S^{1\times n}$.  Denote by $I_r$ ($I$ for short)  the $r\times r$ identity matrix. For $A\in R^{m\times n}$ and  $B\in R^{n\times m}$,
if $AB=I_m$, we call that $A$ has a {\em right inverse} and $B$ is a right inverse of $A$. Similarly, if $AB=I_m$, than  $B$
 has a {\em left inverse} and $A$ is a left inverse of $B$.   The cardinality of a set $X$ is denoted by $|X|$.

A $1\times n$ matrix over $R$ is called  an $n$-dimensional {\em vector} over $R$.
For $\alpha_i\in R^n$, $i=1,\ldots, k$, the vector subset $\{\alpha_1,\ldots, \alpha_k\}$ is called {\em unimodular} if the matrix
$\scriptsize \left(
     \begin{array}{c}
       \alpha_1 \\
       \vdots \\
       \alpha_k \\
     \end{array}
   \right)$ has a right inverse.
Let $V\subseteq R^n$ be a {\em linear subset} (i.e., an $R$-module). A {\em largest unimodular vector set} of $V$ is a unimodular vector subset of $V$
which has  maximum number of vectors. The {\em dimension} of $V$, denoted by ${\rm dim}(V)$, is  the number of vectors in a largest unimodular vector set of $V$.

A linear subset $X$ of $R^n$ is called a {\em $k$-dimensional vector subspace} ({\em $k$-subspace} or {\em subspace} for short) of $R^n$,
if $X$ has a basis  $\{\alpha_1, \alpha_2, \ldots, \alpha_k\}$ being unimodular.
Every basis of a  subspace can be extended to a basis of $R^n$ (cf. \cite[Corollary I.5]{mcdonald2}). We define the $0$-subspace to be $\{0\}$.

Let  $\mathbb{Z}_{p^s}$ denote the {\em residue class ring} of integers modulo $p^s$, where $p$ is a prime and  $s$ is a positive integer.
 The  $\mathbb{Z}_{p^s}$ is a {\em Galois ring}, a commutative {\em local ring},  a finite {\em principal ideal ring} (cf. \cite{mcdonald1, GaloisRing}),
 and a {\em Hermite ring} (according to Cohn's definition \cite{Freeideal,cohn2}). The principal ideal $(p)$ is the unique maximal ideal  of $\mathbb{Z}_{p^s}$,
 and denoted by $J_{p^s}$. The $J_{p^s}$ is also the {\em Jacbson radical} \cite{Bini,Freeideal} of $\mathbb{Z}_{p^s}$.
When $s=1$,  $\mathbb{Z}_{p}$ is a finite field.
For any $x\in \mathbb{Z}_{p^s}$, $x$ is invertible (i.e., a unit) if and only if  $x\notin J_{p^s}$.
We have (cf. \cite{mcdonald1,GaloisRing}) that
\begin{equation}\label{d13fhmm8684}
\left|\mathbb{Z}_{p^s}\right|=p^s, \ \  \left|\mathbb{Z}_{p^s}^*\right|=(p-1)p^{s-1}, \ \ \left|J_{p^s}\right|=p^{s-1}.
\end{equation}

The residue class ring $\mathbb{Z}_{p^s}$ plays  an important role in mathematics and information science. However, since there are zero divisors,
the properties of the vector spaces over $\mathbb{Z}_{p^s}$ are essentially different from the vector spaces over a finite field.
For example, the intersection of two vector subspaces over $\mathbb{Z}_{p^s}$ may not be a vector subspace (cf. Section 3 below).
This brings difficulty to study of geometry and graph theory on $\mathbb{Z}_{p^s}$.  By using matrix representation and the inner rank \cite{cohn2} of a matrix,
we will discuss the basic properties of vector spaces over $\mathbb{Z}_{p^s}$. For example, intersection,  join, dimension formula and
dual subspaces of  $\mathbb{Z}_{p^s}^n$.

Recently, some scholars studied  some graphs (for instance, symplectic graphs and bilinear forms graph) on $\mathbb{Z}_{p^s}$ or a finite commutative ring
(cf. \cite{Gu,Huang-2017,F.Li3,F.Li1,Meemark2}). Note that the Grassmann graph over a finite field plays an important role in geometry \cite{M.Pankov02,W1, Wbook},
graph theory \cite{Brouwera2} and coding theory \cite{Etzion, Etzion02, M.Pankov02}. Thus, the study of Grassmann graph over
$\mathbb{Z}_{p^s}$ has good significance for geometry, combinatorics  and coding theory.

We define the  Grassmann graph over  $\mathbb{Z}_{p^s}$ as follows.
Suppose $m,n$ are integers with $1\leq m< n$. The {\em Grassmann graph} over $\mathbb{Z}_{p^s}$, denoted by  $G_{p^s}(n,m)$, has vertex set the set all $m$-subspaces
of $\mathbb{Z}_{p^s}^n$, and two vertices are adjacent if their intersection is a linear subset of dimension $m-1$.  Based on our results on vector subspaces of $\mathbb{Z}_{p^s}^n$,
we will study  Grassmann graph over $\mathbb{Z}_{p^s}$ and its automorphisms.

The  paper is organized as follows. In Section 2, we recall the basic properties of matrices over $\mathbb{Z}_{p^s}$  and the bilinear forms graphs over $\mathbb{Z}_{p^s}$.
In Section 3,  we will discuss vector subspaces of $\mathbb{Z}_{p^s}^n$. For example,  intersection, join and dimension formula of
two vector subspaces of $\mathbb{Z}_{p^s}^n$;  dual subspaces and arithmetic distances of subspaces.
In Section 5, we show that the Grassmann graph $G_{p^s}(n,m)$ is a connected vertex-transitive graph, and
determine its valency,  clique number and maximum cliques. In Section 6, we characterize the automorphisms of $G_{p^s}(n,m)$.

\section{Matrices and bilinear forms graphs over $\mathbb{Z}_{p^s}$}
 \setcounter{equation}{0}

\ \ \ \ \  In this section,  we recall some basic properties of matrices over $\mathbb{Z}_{p^s}$. For instance, matrix factorization,
the inner rank and McCoy rank of matrix over $\mathbb{Z}_{p^s}$.  We also introduce the bilinear forms graph over $\mathbb{Z}_{p^s}$.

\subsection{Matrices}

\ \ \ \ \
Let  $GL_n(R)$ be the set of $n\times n$ invertible matrices over $R$.
For $A\in R^{m\times n}$, let $^tA$ denote the transpose matrix of a matrix $A$ and ${\rm det}(X)$ the determinant \cite{brown} of a
square  matrix $X$ over $R$.
Let  $0_{m,n}$ ($0$ for short) be the $m\times n$ zero matrix.
Denote by  $E^{m\times n}_{ij}$ ($E_{ij}$ for short)  the $m\times n$ matrix whose $(i,j)$-entry is $1$ and all other entries are $0$'s.
Let ${\rm diag}(A_1,\ldots, A_k)$ denote a block diagonal matrix where $A_i$ is an $m_i\times n_i$ matrix.

For $0\neq A\in R^{m\times n}$, by Cohn's definition \cite{cohn2}, the {\em inner rank} $\rho(A)$ of $A$,  is the least integer $r$ such that
\begin{equation}\label{innerrank1}
\mbox{$A=BC$ \  where  $B\in R^{m\times r}$ and $C\in R^{r\times n}$.}
\end{equation}
 Let $\rho(0)=0$. Any factorization such that as (\ref{innerrank1}) where $r=\rho(A)$ is called a {\em minimal factorization} of $A$.
For $A\in R^{m\times n}$, it is clear that $\rho(A)\leq \min\left\{m,n\right\}$ and $\rho(A)=0$ if and only if $A=0$.
When $R$ is a field, we have $\rho(A)={\rm rank}(A)$, where ${\rm rank}(A)$ is the usual rank of  matrix $A$ over a field.

For matrices over $R$, the followings hold (cf. \cite[Section 5.4]{cohn2,Freeideal}):
\begin{equation}\label{rank-06}
\mbox{$\rho(A)=\rho(PAQ)$ \ where $P$ and $Q$ are invertible matrices over $R$;}
\end{equation}
\begin{equation}\label{rank-02}
\rho(AB)\leq \min\left\{ \rho(A), \rho(B) \right\};
\end{equation}
\begin{equation}\label{rank-04}
\rho(A_1,A_2)\geq\max\left\{\rho(A_1), \rho(A_2) \right\};
\end{equation}
\begin{equation}\label{rank-05}
\rho\left(
            \begin{array}{cc}
              A & * \\
              * & * \\
            \end{array}
          \right)\geq \rho(A).
\end{equation}

Denote by $I_k(A)$ the ideal in $R$ generated by all $k\times k$ minors of $A$, $k=1,\ldots, {\rm min}\{m,n\}$.
Let ${\rm Ann}_R(I_k(A))=\left\{x\in R: xI_k(A)=0\right\}$ denote the annihilator of $I_k(A)$. The {\em McCoy rank} of $A$, denoted by ${\rm rk}(A)$, is the following integer:
$${\rm rk}(A)={\rm max}\left\{k: {\rm Ann}_R(I_k(A))=(0)\right\}.$$
 We have that ${\rm rk}(A)={\rm rk}(^tA)$;
${\rm rk}(A)={\rm rk}(PAQ)$ where $P$ and $Q$ are invertible matrices of the appropriate sizes; and ${\rm rk}(A)=0$ if and only if ${\rm Ann}_R(I_1(A))\neq (0)$ (cf. \cite{brown, mcdonald4}).

 Since $\mathbb{Z}_{p^s}$ is an {\em Hermite ring} (according to Cohn's definition \cite{Freeideal}),
every matrix over $\mathbb{Z}_{p^s}$ which has  right inverse is completable (i.e., it can be  extended to an invertible matrix over $\mathbb{Z}_{p^s}$).
Since $\mathbb{Z}_{p^s}$ is a commutative local ring, a matrix $\scriptsize \left(
     \begin{array}{c}
       \alpha_1 \\
       \vdots \\
       \alpha_m \\
     \end{array}
   \right)\in\mathbb{Z}_{p^s}^{m\times n}$  has a right inverse if and only if $\alpha_1, \ldots, \alpha_m$ are $m$ linearly independent vectors over $\mathbb{Z}_{p^s}$.


\begin{lem}\label{canonical}
{\rm (cf. \cite[Chap. II]{Newman}), or \cite[p.327]{mcdonald1}, or \cite[Lemma 2.3]{Huang-2017})} \
Let $R=\mathbb{Z}_{p^s}$, and let $A\in R^{m\times n}$ be a non-zero matrix.
Then there are $P\in GL_m(R)$ and $Q\in GL_n(R)$ such that
\begin{equation}\label{form1}
A=P\left(\begin{array}{ccccc}
          I_r & & & & \\
              & p^{k_1}& & & \\
              & & \ddots & & \\
              & & & p^{k_t}&  \\
              & & & & 0 \\
        \end{array}\right)Q,
\end{equation}
where  $1\leq k_1\leq \cdots\leq k_t\leq {\rm max}\{s-1,1\}$. Moreover, the parameters $(r, t, k_1,\ldots, k_t)$ are uniquely determined by $A$.
In (\ref{form1}), $I_r$ or ${\rm diag}\left(p^{k_1},\ldots,p^{k_t}\right)$ may be absent.
\end{lem}

\begin{lem}\label{innerrank02} {\rm(see \cite[Lemmas 2.4 and 2.7]{Huang-2017})} \
Let $0\neq A\in \mathbb{Z}_{p^s}^{m\times n}$ ($s\geq 2$) be of the form (\ref{form1}). Then $r+t$ is  the inner rank of $A$,
and $r$ is  the McCoy rank of $A$.
\end{lem}

Let $A\in \mathbb{Z}_{p^s}^{m\times n}$ ($s\geq 2$). By Lemma \ref{innerrank02}, ${\rm rk}(A)\leq \rho(A)$, and
${\rm rk}(A)=0$ if and only if $A\in  J_{p^s}^{m\times n}$.
For matrices $A,B,C$ over $\mathbb{Z}_{p^s}$, By Lemmas \ref{canonical} and \ref{innerrank02}, it is easy to see that (cf. \cite{Huang-2017})
\begin{equation}\label{rank-0008}
\mbox{$\rho\left(
            \begin{array}{cc}
              A & 0 \\
              0 & B\\
            \end{array}
          \right)=\rho(A)+\rho(B)$;}
\end{equation}
\begin{equation}\label{rank-7}
 \rho(A+C)\leq \rho(A)+\rho(C).
\end{equation}

\begin{lem}\label{gl}{\rm (see  \cite[Corollary 2.21]{brown})} \ Let $A\in \mathbb{Z}_{p^s}^{n\times n}$.
 Then $A\in GL_n(\mathbb{Z}_{p^s})$ if and only if $\det(A)\in \mathbb{Z}_{p^s}^*$.
\end{lem}

Note that if $a\in\mathbb{Z}_{p^s}^*$ and $b\in J_{p^s}$, then $a\pm b\in\mathbb{Z}_{p^s}^*$.
Lemma \ref{canonical}  implies  the following lemma.

\begin{lem}\label{Mc-rank-1}{\rm (cf. \cite[Lemma 2.10]{Huang-2017})} \
Let $A\in \mathbb{Z}_{p^s}^{m\times n}$ where $s\geq 2$ and $n\geq m$. Then ${\rm rk}(A)={\rm rk}(A\pm B)$ for all $B\in J_{p^s}^{m\times n}$.
Moreover, $A$ has a right inverse if and only if ${\rm rk}(A)=m$.
\end{lem}

Let $T_p=\{0,1,\ldots, p-1\}\subseteq\mathbb{Z}_{p^s}$.
For two distinct elements $a,b\in T_p$,  we have $a-b\in \mathbb{Z}_{p^s}^*$.
Without loss of generality, we may assume that $T_p=\mathbb{Z}_p$ in our discussion.

\begin{lem}\label{uniqueuk1}{\rm (see \cite[Proposition 6.2.2]{Bini} or \cite[p.328]{mcdonald1})} \ Every non-zero element $x$ in $\mathbb{Z}_{p^s}$ can be written as
$x=up^t$ where $u$ is a unit and  $0\leq t\leq s-1$. Moreover, the integer $t$ is unique and $u$ is unique modulo the ideal $(p^{s-t})$ of $\mathbb{Z}_{p^s}$.
\end{lem}

\begin{lem}\label{uniqueuk2}{\rm (cf. \cite[p.328]{mcdonald1})} \  Every non-zero element $x$ in $\mathbb{Z}_{p^s}$ can be written uniquely as
$$\mbox{$x=t_0+t_1p+\cdots +t_{s-1}p^{s-1}$, \ where $t_i\in T_p$,  $i=0,1,\ldots, s-1$.}$$
\end{lem}

By Lemma \ref{uniqueuk2},
every matrix $X\in\mathbb{Z}_{p^s}^{m\times n}$ can be written uniquely as
\begin{equation}\label{bcweu7754332t}
\mbox{$X=X_0+X_1p+\cdots +X_{s-1}p^{s-1}$, \ where $X_i\in T_p^{m\times n}$, $i=0,\ldots, s-1$.}
\end{equation}

Note that every matrix in $T_p^{m\times n}$ can be seen as a matrix in $\mathbb{Z}_p^{m\times n}$.
We define the {\em natural surjection}
\begin{equation}\label{naturalsurjection}
\pi: \mathbb{Z}_{p^s}^{m\times n}\rightarrow \mathbb{Z}_{p}^{m\times n}
\end{equation}
by $\pi(X)=X_0$ for all $X\in\mathbb{Z}_{p^s}^{m\times n}$ of the form (\ref{bcweu7754332t}). Clearly, $\pi(A)=A$ if $A\in\mathbb{Z}_p^{m\times n}$.
For  $X, Y\in\mathbb{Z}_{p^s}^{m\times n}$ and $Q\in\mathbb{Z}_{p^s}^{n\times k}$, We have
\begin{equation}\label{pipipipi003a}
 \pi(X+Y)=\pi(X)+\pi(Y),
 \end{equation}
\begin{equation}\label{pipipipi003}
 \pi(XQ)=\pi(X)\pi(Q),
\end{equation}
\begin{equation}\label{pipipipi004}
 \pi(^tX)= \,^t(\pi(X)).
\end{equation}

\begin{thm}\label{pipipipi005} \ Let  $X\in\mathbb{Z}_{p^s}^{m\times n}$ $(n\geq m$). Then the following hold:
\begin{itemize}
\item[{\rm (i)}]   ${\rm rk}(X)={\rm rank}(\pi(X))$.

\item[{\rm (ii)}]   $X$ has a right inverse if and only if $\pi(X)$ has a right inverse.

\item[{\rm (iii)}]   When $n=m$, $X$ is invertible  if and only if $\pi(X)$ is invertible. Moreover, if $X$ is invertible, then
\begin{equation}\label{pipipipi006}
 \pi(X^{-1})= (\pi(X))^{-1}.
 \end{equation}
 \end{itemize}
\end{thm}
\proof
Without loss of generality, we assume that $s\geq 2$.

(i). \ Let ${\rm rk}(X)=r$. By Lemma \ref{canonical}, there are matrices $P\in GL_m(\mathbb{Z}_{p^s})$ and $Q\in GL_n(\mathbb{Z}_{p^s})$ such that
$$X=P{\rm diag}\left(I_r, D,0_{m-r-t, n-r-t} \right)Q,$$
 where $D={\rm diag}\left(p^{k_1},\ldots, p^{k_t}\right)$,  $1\leq k_1\leq \cdots\leq k_t\leq s-1$.
Using (\ref{pipipipi003}), we have that $\pi(P)$ and $\pi(Q)$ are invertible, and
$\pi(X)=\pi(P){\rm diag}\left(I_r, 0_t, 0_{m-r-t, n-r-t} \right)\pi(Q)$.
Thus ${\rm rank}(\pi(X))=r={\rm rk}(X)$.

(ii). \ If $X$ has a right inverse, then (\ref{pipipipi003}) implies that $\pi(X)$ has a right inverse. Conversely, suppose  $\pi(X)$ has a right inverse.
Using (i), we get ${\rm rank}(\pi(X))=m={\rm rk}(X)$. If follows from Lemma \ref{Mc-rank-1} that $X$ has a right inverse.

(iii). \ Suppose that $n=m$. Similar to the proof of (ii), $X$ is invertible  if and only if $\pi(X)$ is invertible.
Now, let $X=X_0+X_1p+\cdots +X_{s-1}p^{s-1}\in\mathbb{Z}_{p^s}^{n\times n}$ be invertible, where $X_i\in\mathbb{Z}_p^{n\times n}$, $i=0,\ldots, s-1$.
Write $X^{-1}=Y_0+Y_1p+\cdots +Y_{s-1}p^{s-1}\in\mathbb{Z}_{p^s}^{n\times n}$, where $Y_i\in\mathbb{Z}_p^{n\times n}$, $i=0,\ldots, s-1$.
Then $X_0Y_0=I_n$, and hence (\ref{pipipipi006}) holds.
$\qed$

\begin{lem}\label{ds424rrwrww}{\rm (see \cite[Lemma 4.2]{Huang-2017})} \  If $A\in GL_n(\mathbb{Z}_{p^{s-1}})$ where $s\geq 2$, then $A\in GL_n(\mathbb{Z}_{p^{s}})$.
 \end{lem}

\begin{lem}\label{64gdgd7wrww}{\rm (see \cite[Lemma 4.3]{Huang-2017})} \  If  $A\in\mathbb{Z}_{p^{s-1}}^{m\times n}$ where $s\geq 2$, then both $A$ and $Ap$
can be viewed as matrices in $\mathbb{Z}_{p^{s}}^{n\times n}$ with the same inner rank.
\end{lem}

\subsection{Bilinear forms graphs}

\ \ \ \ \ All  graphs in this paper are {\em simple} \cite{Godsil} and finite. Suppose  $G$ is a graph.
 Denote the {\em distance} between vertices $x$ and $y$ in $G$ by $d_G(x, y)$ ($d(x,y)$ for short).
 Let $V(G)$ be the vertex set of $G$. For $x,y\in V(G)$, we write $x \sim y$ if vertices $x$ and $y$  are adjacent.

Recall that a {\em clique} of $G$ is a complete subgraph of $G$. A clique $C$ is  maximal if there is no clique of $G$ which properly contains $C$ as a subset.
A {\em maximum clique} of $G$ is a clique of $G$ which has maximum cardinality.
The  {\em clique number} $\omega(G)$ of $G$  is the number of vertices in a  maximum clique.
For convenience, we regard that a maximal clique and its vertex set are the same.

An {\em independent set} of $G$ is a subset of vertices such that no two vertices are adjacent. A {\em largest independent set}
of $G$ is an independent set of maximum cardinality. The {\em independence number} $\alpha(G)$ is  the number of vertices in a largest independent set of $G$.

As a natural extension of the bilinear forms graph over a finite field,  the bilinear forms graph over  $\mathbb{Z}_{p^s}$ was defined as follows (cf. \cite{Huang-2017}).
The {\em bilinear forms graph} over $\mathbb{Z}_{p^s}$, denoted by  $\Gamma(\mathbb{Z}_{p^s}^{m\times n})$,
has the vertex set $\mathbb{Z}_{p^s}^{m\times n}$ where $m,n\geq 2$, and two vertices $A$ and $B$ are adjacent if $\rho(A-B)=1$.

In $\mathbb{Z}_{p^{s}}^{m\times n}$, let
$$\mbox{$\mathcal{M}_1=\left\{\sum\limits_{j=1}^nx_jE_{1j}: x_j\in \mathbb{Z}_{p^s} \right\}$, \ \ $\mathcal{N}_1=\left\{\sum\limits_{i=1}^mx_iE_{i1}: x_i\in \mathbb{Z}_{p^s} \right\}$.}$$
 A maximal clique $\mathcal{C}$ of $\mathbb{Z}_{p^{s}}^{m\times n}$ is called of {\em type one} (resp. {\em type two}),
if $\mathcal{C}$ is of the form $\mathcal{C}=P\mathcal{M}_1+A$ (resp. $\mathcal{C}=\mathcal{N}_1Q+A$),
where $P$ and $Q$ are fixed invertible matrices over $\mathbb{Z}_{p^s}$ and $A\in\mathbb{Z}_{p^s}^{m\times n}$ is fixed.

\begin{lem}\label{gdf5435dt000}{\rm(see \cite[Lemma 3.4]{Huang-2017})} \
Let $\mathcal{C}$ be a maximal clique of $\Gamma(\mathbb{Z}_{p^s}^{m\times n})$ ($s\geq 2$) containing the vertex $0$.
Then $\mathcal{C}\nsubseteq J_{p^s}^{m\times n}$.
\end{lem}

\begin{lem}\label{gdf5435dt}{\rm(see \cite[Theorem 3.6]{Huang-2017})} \
 Let $\Gamma=\Gamma(\mathbb{Z}_{p^s}^{m\times n})$ and  $k={\rm max}\{m,n\}$. Then
\begin{equation}\label{rwe4gdg5353gd}
\omega\left(\Gamma\right)=p^{sk}.
\end{equation}
Moreover, if $n>m$,  then every  maximum clique of $\Gamma$ is of type one. If $m=n$,  then every  maximum clique of $\Gamma$ is of type one or type two.
\end{lem}

By \cite[Theorem 4.4]{Huang-2017}, if $2\leq m\leq n$, then
\begin{equation}\label{645rerta12bnn}
\mbox{$\alpha \left(\Gamma(\mathbb{Z}_{p^s}^{m\times n})\right)=p^{sn(m-1)}$.}
\end{equation}


\section{Subspaces of $\mathbb{Z}_{p^s}^n$}
 \setcounter{equation}{0}

\ \ \ \ \ In this section, we study vector subspaces of $\mathbb{Z}_{p^s}^n$. For example,  intersection, join and dimension formula of
two  subspaces of $\mathbb{Z}_{p^s}^n$. These results have many applications in geometry and combinatorics.

In  general, a linear subset $X$ in $R^n$ may not be a subspace, although it has a dimension.
For example, let $R=\mathbb{Z}_{p^s}$, and let $X\subseteq J_{p^s}^n$ be a linear subset with $X\neq \{0\}$. Then ${\rm dim}(X)=0$ because $p^{s-1}x=0$ for all $x\in X$.
 Thus, $X$ is not a subspace of $R^n$. On the other hand, if $X\subseteq \mathbb{Z}_{p^s}^n$ is a linear subset with ${\rm dim}(X)=0$, we cannot imply $X=\{0\}$ (i.e., a $0$-subspace).

Let $X=[\alpha_1,\ldots,\alpha_k]$ be a $k$-subspace of $\mathbb{Z}_{p^s}^n$. Then $X$ has a {\em matrix representation}
$\scriptsize \left(
     \begin{array}{c}
       \alpha_1 \\
       \vdots \\
       \alpha_k \\
     \end{array}
   \right)\in \mathbb{Z}_{p^s}^{k\times n}$  (cf. \cite{L.P.Huang2009-3, Huang-book}).
For simpleness, the matrix representation of  $X$ is also written as $X$.
If $X$ is a matrix representation of the subspace $X$, then for any invertible matrix $P\in GL_k(\mathbb{Z}_{p^s})$,
$PX$ is also a matrix representation of $X$. Thus, the matrix representation is not unique. However,
the subspace $X$ has a unique matrix representation which is the {\em row-reduced echelon form} $X=(I_k, B)Q$, where $Q$ is a permutation matrix.

\subsection{Join and dimensional formula of subspaces}

\ \ \ \ \
In $\mathbb{Z}_{p^s}^n$, a {\em join} of two subspaces $A$ and $B$ is a minimum dimensional subspace containing  $A$ and $B$.
 In general, when $s\geq 2$ the join  is not unique.
For example, suppose that $s\geq 3$, $\alpha=(1,p,0)$, $\beta=(1,p^2,0)$, $\gamma=(0,1,p^{s-1})$,
$e_1=(1,0,0)$ and $e_2=(0,1,0)$. Then $[e_1,e_2]$ and $[e_1, \gamma]$ are two distinct joins of $[\alpha]$ and $[\beta]$.
Denoted by $A\vee B$ the set of all joins of subspaces $A$ and $B$ in $\mathbb{Z}_{p^s}^n$ with the same minimum dimension $\dim(A\vee B)$.
We have $A\vee B=B\vee A$. When $B\subseteq A$ are two subspaces, we have $A\vee B=\{A\}$.

Let $X,Y$ be two subspaces of $\mathbb{Z}_{p^s}^n$.  In general, $X\cap Y$ is a linear subset but it may not be a subspace of $\mathbb{Z}_{p^s}^n$.
For example, let $\scriptsize X=\left(\begin{array}{cccc}
                        1 & 0 & 0 & 0 \\
                        0 & 1 & 0 & 0 \\
                      \end{array}
                    \right)$,
$\scriptsize Y=\left(\begin{array}{cccc}
                        1 & 0 & 0 & 0 \\
                        0 & 1 & p^{s-1} & 0 \\
                      \end{array}
                    \right)$ ($s\geq 2$) are  $2$-subspaces of $\mathbb{Z}_{p^s}^4$.
Then $(1,p^i, 0,0)\in X\cap Y$, $i=0,1,  \ldots, s-1$.   It follows that $X\cap Y$ is not a free  $\mathbb{Z}_{p^s}$-module, and hence  $X\cap Y$ is not a subspace of $\mathbb{Z}_{p^s}^4$.

\begin{thm}\label{hgf5tedrdg3} \ Let $n>k\geq m \geq 1$.  Suppose $A$ and $B$ are $k$-subspace and $m$-subspace of $\mathbb{Z}_{p^s}^n$, respectively, and $B\nsubseteq A$.
Then there is $U\in GL_n(\mathbb{Z}_{p^s})$ such that
\begin{equation}\label{HGFYR234LJK0r}
\mbox{$A=(0, I_k)U$ \ and \ $B=(D, I_m)U$,}
\end{equation}
 where $D={\rm diag}(p^{i_1}, \ldots, p^{i_r}, 0_{m-r, n-m-r})$,  $1\leq r\leq {\rm min}\{m, n-k\}$ and $0\leq i_1\leq \ldots \leq i_r\leq {\rm max}\{s-1, 1\}$.
\end{thm}
\proof
Without loss of generality, we assume that $s\geq 2$. By the matrix representation of subspace and Lemma \ref{canonical},
there is $U_1\in GL_n(\mathbb{Z}_{p^s})$ such that $A=(0, I_k)U_1$. Without loss of generality, we assume $U_1=I_n$ and hence $A=(0, I_k)$. Write $B=(B_1,B'_2)$ where
$B_1\in\mathbb{Z}_{p^s}^{m\times (n-k)}$ and $B'_2\in\mathbb{Z}_{p^s}^{m\times k}$. When $k>m$, by elementary operations of matrix, we may assume with no loss of generality that
$B_2'=(0_{m,k-m},B_2)$ where $B_2\in\mathbb{Z}_{p^s}^{m\times m}$. Thus $B=(B_1, 0_{m,k-m}, B_2)$.
Using  Lemma \ref{canonical},  we may assume with no loss of generality that
$$B_1={\rm diag}(p^{i_1}, \ldots, p^{i_r}, 0_{m-r, n-k-r})$$
 where $1\leq r\leq {\rm min}\{m, n-k\}$ and $0\leq i_1\leq \ldots \leq i_r\leq s-1$.

{\bf Case 1}. \ $B_2$ is invertible.

Let $U_2={\rm diag}(I_{n-m}, B_2^{-1})$. Then $A=(0, I_k)=(0, I_k)U_2$ and $B=(B_1, 0_{m,k-m}, I_m)U_2$.
Thus (\ref{HGFYR234LJK0r}) holds.

{\bf Case 2}. \ $B_2$ is not invertible.

 Without loss of generality, we let
$\mbox{${\rm diag}(p^{i_1}, \ldots, p^{i_{r}})={\rm diag}(I_t, D_1)$}$,
where
$D_1={\rm diag}(p^{i_{t+1}}, \ldots, p^{i_{r}})$, $1\leq i_{t+1}\leq \ldots \leq i_r\leq s-1$.
 Write $\scriptsize B_2=\left(\begin{array}{c}
B_{21} \\
 B_{22} \\
  \end{array}
  \right)\in \mathbb{Z}_{p^s}^{m\times m}$
where $B_{21}\in\mathbb{Z}_{p^s}^{r\times m}$ and $B_{22}\in\mathbb{Z}_{p^s}^{(m-r)\times m}$. Then $B_{22}$ has a right inverse.
Without loss of generality, we may assume that $B_{22}=(I_{m-r},0)$.  Then $A=(0,I_k)$ and
$$B=\left(
      \begin{array}{cccc}
        I_t & 0 & 0 &Y_{1}\\
        0 & D_2 & 0& Y_{2} \\
        0 & 0 & I_{m-r} & 0 \\
      \end{array}
    \right),$$
where $D_2$ is an $(r-t)\times (n-m-t)$ diagonal matrix over $J_{p^s}$, $Y_1\in\mathbb{Z}_{p^s}^{t\times r}$ and $Y_2\in\mathbb{Z}_{p^s}^{(r-t)\times r}$.
Since $(0, D_2, 0, Y_2)$ has a right inverse,  it is easy to see that $Y_2$ has a right inverse. Hence
 there exists $C_1\in\mathbb{Z}_{p^s}^{t\times r}$ such that
$\scriptsize \left( \begin{array}{c}
 C_1+Y_1 \\
 Y_2 \\
   \end{array}
  \right)$  is invertible. Put
  $$U_3=\left(
      \begin{array}{cccc}
        I_t & 0 & 0 & C_1\\
        0 & I_{n-m-t} & 0& 0 \\
        0 & 0 & I_{m-r} & 0 \\
        0 & 0 & 0 & I_{r} \\
      \end{array}
    \right)\in GL_n(\mathbb{Z}_{p^s}).$$
 Then $A=AU_3$ and
 $$\small BU_3=\left(
      \begin{array}{cccc}
        I_t & 0 & 0 &Y_{1}\\
        0 & D_2 & 0& Y_{2} \\
        0 & 0 & I_{m-r} & 0 \\
      \end{array}
    \right)\left(
      \begin{array}{cccc}
        I_t & 0 & 0 & C_1\\
        0 & I_{n-m-t} & 0& 0 \\
        0 & 0 & I_{m-r} & 0 \\
        0 & 0 & 0 & I_{r} \\
      \end{array}
    \right)=\left(
      \begin{array}{cccc}
        I_t & 0 & 0 &C_1+Y_{1}\\
        0 & D_2 & 0& Y_{2} \\
        0 & 0 & I_{m-r} & 0 \\
      \end{array}
    \right).
 $$
Since  $\scriptsize B_2'':= \left(\begin{array}{cc}
         0 &C_1+Y_{1}\\
          0& Y_{2} \\
         I_{m-r} & 0 \\
      \end{array}
    \right)$ is invertible,   Case 1 implies that (\ref{HGFYR234LJK0r}) holds.
$\qed$

\begin{thm}\label{dimension-formula01}{\rm(dimensional formula)} \ Let $A$ and $B$ be two subspaces of $\mathbb{Z}_{p^s}^n$.
Then
\begin{equation}\label{dimension-formula02}
{\rm dim}(A\vee B)={\rm dim}(A)+{\rm dim}(B)-{\rm dim}(A\cap B)= \rho\left(
\begin{array}{c}
A \\
B \\
\end{array}
\right).
\end{equation}
\end{thm}
\proof
Let ${\rm dim}(A)=k$ and ${\rm dim}(B)=m$. Without loss of generality, we assume that $s\geq 2$,  $n>k\geq m\geq 1$ and $B\nsubseteq A$.
By Theorem \ref{hgf5tedrdg3}, there is $U\in GL_n(\mathbb{Z}_{p^s})$ such that
$A=(0, I_k)U$ \ and \ $B=(D, I_m)U$,
where $D={\rm diag}(p^{i_1}, \ldots, p^{i_r}, 0_{m-r, n-m-r})$,  $1\leq r\leq {\rm min}\{m, n-k\}$ and $0\leq i_1\leq \ldots \leq i_r\leq s-1$.
Put  $D_1={\rm diag}(p^{i_1}, \ldots, p^{i_r})$.
Without loss of generality, we may assume that $U=I_n$ and $m\geq 2$.  Hence  we can  assume further that
$$A=(0,I_k), \ \ \ B=\left(\begin{array}{cccc}
                           D_1 & 0_{r,n-m-r} & I_r & 0 \\
                           0 & 0 & 0 & I_{m-r} \\
                         \end{array}
                       \right).$$
In above matrices, some zero elements of matrices may be absent.

Clearly, the $(m-r)$-subspace $(0,0,0, I_{m-r})$ is contained in $A\cap B$. On the other hand,  the $(k+r)$-subspace
$\small \left( \begin{array}{ccccc}
                           I_r & 0_{r,n-k-r}& 0 \\
                           0  & 0 & I_{k} \\
                         \end{array}
                       \right)$ contains $A$ and $B$.
It follows that
\begin{equation}\label{dimension-formula04}
{\rm dim}(A\cap B)\geq m-r, \ \ \ {\rm dim}(A\vee B)\leq k+r.
\end{equation}

Let $\alpha\in A\cap B$ be an $n$-dimensional vector. Then there are matrices $T_1\in\mathbb{Z}_{p^s}^{1\times k}$ and $T_2=(T_{21},T_{22})\in\mathbb{Z}_{p^s}^{1\times m}$,
 where $T_{21}\in\mathbb{Z}_{p^s}^{1\times r}$ and  $T_{22}\in\mathbb{Z}_{p^s}^{1\times (m-r)}$, such that
$\alpha=T_1A=(0, T_1)$,  $\alpha=T_2B=(T_{21}D_1, 0,T_{21}, T_{22})$.
Hence $T_{21}D_1=0$ and $\alpha=(0, 0,T_{21}, T_{22})$. By  $T_{21}D_1=0$ we get $T_{21}\in J_{p^s}^{1\times r}$. Thus,
$$\alpha=(0,\ldots,0,t_{1}, \ldots, t_{r}, y_1, \ldots, y_{m-r})$$
where $t_1, \ldots,t_r\in J_{p^s}$.
It follows that ${\rm dim}(A\cap B)\leq m-r$. Hence, from (\ref{dimension-formula04}) we obtain
\begin{equation}\label{dimension-formula05}
{\rm dim}(A\cap B)=m-r.
\end{equation}

Assume that $d={\rm dim}(A\vee B)$. Then there exists a $d$-subspace $W\in A\vee B$. Since $B\nsubseteq A$, one has $d>k$.
Thus $\scriptsize W=\left(\begin{array}{c}
                           W_1 \\
                           A  \\
                         \end{array}
                       \right)=\left(\begin{array}{cc}
                           W_{11} & 0 \\
                           0& I_k  \\
                         \end{array}
                       \right)$, where $W_{11}\in\mathbb{Z}_{p^s}^{(d-k)\times (n-k)}$ has a right inverse. Since $B\subset W$, the $r$-subspace
$B_1:=(D_1,0_{r,n-m-r},I_r,0)\subset W$. Thus, there is a matrix $P=(P_1,P_2)\in \mathbb{Z}_{p^s}^{r\times d}$ where $P_1\in \mathbb{Z}_{p^s}^{r\times (d-k)}$ and
$P_2\in\mathbb{Z}_{p^s}^{r\times k}$, such that $B_1=PW$ and hence $(D_1,0_{r,n-m-r},I_r,0)=(P_1W_{11}, P_2)$.
It follows from (\ref{rank-02}) that
$$r=\rho(D_1,0_{r,n-k-r})=\rho(P_1W_{11})\leq \rho(W_{11})\leq d-k.$$
Hence $d\geq k+r$. By (\ref{dimension-formula04}) we obtain
\begin{equation}\label{dimension-formula06}
{\rm dim}(A\vee B)=k+r.
\end{equation}

Clearly, we have $\scriptsize \rho\left(
\begin{array}{c}
A \\
B \\
\end{array}
\right)=\rho\left(
\begin{array}{ccc}
0&0&I_k \\
D_1&0 & 0 \\
0&0 & 0 \\
\end{array}
\right)=k+r$. By (\ref{dimension-formula05}) and (\ref{dimension-formula06}), we get (\ref{dimension-formula02}).
$\qed$

\begin{rem}\label{dimension-formula08} \ Let $A$ and $B$ be two subspaces of $\mathbb{Z}_{p^s}^n$. By (\ref{dimension-formula02}), $A\vee B$ is the set of subspaces
containing $A$ and $B$ with the same dimension  $\scriptsize \rho\left(
\begin{array}{c}
A \\
B \\
\end{array}
\right)$.
 \end{rem}

\subsection{On intersection and join of subspaces}

\begin{thm}\label{fd32nvnafal87990} \  Let $A$ and $B$ be two subspaces of $\mathbb{Z}_{p^s}^n$ and ${\rm dim}(A)+{\rm dim}(B)\leq n$.
Then ${\rm dim}(A\cap B)=0$ if and only if
$\scriptsize\rho\left(
\begin{array}{c}
A \\
B \\
\end{array}
\right)=\rho(A)+\rho(B)$. Moreover, $A\cap B=\{0\}$ if and only if
$\scriptsize {\rm rk}\left(
\begin{array}{c}
A \\
B \\
\end{array}
\right)=\rho(A)+\rho(B)$.
\end{thm}
\proof
Let ${\rm dim}(A)=k=\rho(A)$ and ${\rm dim}(B)=m=\rho(B)$.
Without loss of generality, we assume that $s\geq 2$,  $k\geq m\geq 1$ and $B\nsubseteq A$. Clearly, $\scriptsize {\rm rk}\left(
\begin{array}{c}
A \\
B \\
\end{array}
\right)\leq\rho\left(
\begin{array}{c}
A \\
B \\
\end{array}
\right)\leq k+m$.

Suppose $\scriptsize\rho\left(
\begin{array}{c}
A \\
B \\
\end{array}
\right)<k+m$.
By Theorem \ref{hgf5tedrdg3}, there is $U_1\in GL_n(\mathbb{Z}_{p^s})$ such that $A=(0, I_k)U_1$ and $B=(D, I_m)U_1$,
 where $D={\rm diag}(p^{i_1}, \ldots, p^{i_r}, 0_{m-r, n-m-r})$,  $1\leq r<m$ and $0\leq i_1\leq \ldots \leq i_r\leq s-1$.
 Thus  $\alpha:=(0, \ldots, 0,1)U_1\in A\cap B$. Since $\alpha$ is unimodular,   ${\rm dim}(A\cap B)\geq 1$.
Conversely, if ${\rm dim}(A\cap B)\geq 1$, then it is clear that $\scriptsize\rho\left(
\begin{array}{c}
A \\
B \\
\end{array}
\right)<k+m$. Therefore, ${\rm dim}(A\cap B)\geq 1$ if and only if $\scriptsize\rho\left(
\begin{array}{c}
A \\
B \\
\end{array}
\right)<k+m$. It follows that ${\rm dim}(A\cap B)=0$ if and only if $\scriptsize\rho\left(
\begin{array}{c}
A \\
B \\
\end{array}
\right)=k+m$.

If $\scriptsize {\rm rk}\left(
\begin{array}{c}
A \\
B \\
\end{array}
\right)=k+m$, then  Theorem \ref{hgf5tedrdg3} implies that there is $U_2\in GL_n(\mathbb{Z}_{p^s})$ such that $A=(0, 0,I_k)U_2$ and $B=(I_m,0, I_m)U_2$.
Thus, it is easy to see that  $A\cap B=\{0\}$.

Now, we assume that $A\cap B=\{0\}$. Then ${\rm dim}(A\cap B)=0$, and hence $\scriptsize\rho\left(
\begin{array}{c}
A \\
B \\
\end{array}
\right)=k+m$ because (\ref{dimension-formula02}). Suppose $\scriptsize {\rm rk}\left(
\begin{array}{c}
A \\
B \\
\end{array}
\right)=k+t<k+m$.
By Theorem \ref{hgf5tedrdg3}, there is $U_3\in GL_n(\mathbb{Z}_{p^s})$ such that $A=(0,0, I_k)U_3$ and $B=(D_1, 0, I_m)U_3$, where
$D_1={\rm diag}(I_t, p^{i_{t+1}}, \ldots, p^{i_m})$, $1\leq i_{t+1}\leq \ldots \leq i_m\leq s-1$. Let $\beta=(0, \ldots, 0,p^{s-i_m})U_3$.
Then $\beta\neq 0$. By $\beta=p^{s-i_m}(0, \ldots, 0,p^{i_m}, 0,\ldots, 0,1)U_3=p^{s-i_m}(0, \ldots, 0, 1)U_3$,
we get $\beta\in A\cap B$, a contradiction. Hence we must have
$\scriptsize {\rm rk}\left(
\begin{array}{c}
A \\
B \\
\end{array}
\right)=k+m$. Thus, we have that $A\cap B=\{0\}$ if and only if
$\scriptsize {\rm rk}\left(
\begin{array}{c}
A \\
B \\
\end{array}
\right)=k+m$.
$\qed$

\begin{thm}\label{08ujivxwfal879u} \  Let $A$ and $B$ be two subspaces  of $\mathbb{Z}_{p^s}^n$. Then:
\begin{itemize}
\item[{\rm (i)}] $A\cap B$ is a fixed subspace of dimension ${\rm dim}(A\cap B)$ if and only if
$\scriptsize \rho\left(
\begin{array}{c}
A \\
B \\
\end{array}
\right)={\rm rk}\left(
\begin{array}{c}
A \\
B \\
\end{array}
\right)$.

\item[{\rm (ii)}] $A\vee B$ is a fixed subspace of dimension ${\rm dim}(A\vee B)$ if and only if
$\scriptsize \rho\left(
\begin{array}{c}
A \\
B \\
\end{array}
\right)={\rm rk}\left(
\begin{array}{c}
A \\
B \\
\end{array}
\right)$ or ${\rm dim}(A\vee B)=n$.
\end{itemize}
\end{thm}
\proof
 Put ${\rm dim}(A)=k$ and ${\rm dim}(B)=m$.
Without loss of generality, we assume that $s\geq 2$, $n>k\geq m \geq 1$ and $B\nsubseteq A$.
By Theorem \ref{hgf5tedrdg3}, there is $U\in GL_n(\mathbb{Z}_{p^s})$ such that
$A=(0, I_k)U$  and  $B=(D, I_m)U$,
where $D={\rm diag}(p^{i_1}, \ldots, p^{i_r}, 0_{m-r, n-m-r})$,  $1\leq r\leq {\rm min}\{m, n-k\}$ and $0\leq i_1\leq \ldots \leq i_r\leq s-1$.
Thus $\scriptsize\rho\left(
\begin{array}{c}
A \\
B \\
\end{array}
\right)=k+r$ and hence $n\geq k+r$. Without loss of generality, we may assume $U=I_n$.
By Theorem \ref{dimension-formula01}, we have ${\rm dim}(A\vee B)=k+r$ and ${\rm dim}(A\cap B)=m-r$.

(i). \ Suppose $\scriptsize\rho\left(
\begin{array}{c}
A \\
B \\
\end{array}
\right)={\rm rk}\left(
\begin{array}{c}
A \\
B \\
\end{array}
\right)$. Then $D={\rm diag}(I_r, 0_{m-r, n-m-r})$, $\scriptsize B=\left(
                                                            \begin{array}{cccc}
                                                              I_r & 0 &  I_r &0\\
                                                               0& 0 & 0&I_{m-r} \\
                                                            \end{array}
                                                          \right)$ and
$A=(0,I_k)$. Thus $A\cap B$ contains an $(m-r)$-subspace $A_1:=(0,0,0,I_{m-r})$,
For any vector $x\in A\cap B$, it is easy to verify $x\in A_1$. Therefore, $A\cap B=A_1$ is a fixed $(m-r)$-subspace.

Conversely, suppose $A\cap B$ is a fixed $(m-r)$-subspace. We show $\scriptsize\rho\left(
\begin{array}{c}
A \\
B \\
\end{array}
\right)={\rm rk}\left(
\begin{array}{c}
A \\
B \\
\end{array}
\right)$. Otherwise, $\scriptsize\rho\left(
\begin{array}{c}
A \\
B \\
\end{array}
\right)>{\rm rk}\left(
\begin{array}{c}
A \\
B \\
\end{array}
\right)$, and hence $0\neq p^{i_r}\in J_{p^s}$ in the matrix $D$. By row elementary operations of matrix,
it is easy to see that the $(m-r)$-subspace
$\small A_2:=\left(
                                                            \begin{array}{cccc}
                                                              0 & p^{s-i_r} & 1 &0\\
                                                               0& 0 & 0&I_{m-r-1} \\
                                                            \end{array}
                                                          \right)\subseteq A\cap B$.
Since the $(m-r)$-subspace $A_1=(0,0,0,I_{m-r})\subseteq A\cap B$ and $A_1\neq A_2$, we get a contradiction.

(ii). \  Suppose $\scriptsize\rho\left(
\begin{array}{c}
A \\
B \\
\end{array}
\right)={\rm rk}\left(
\begin{array}{c}
A \\
B \\
\end{array}
\right)$. Similarly,  $\scriptsize B=\left(
                                                            \begin{array}{cccc}
                                                              I_r & 0 &  I_r &0\\
                                                               0& 0 & 0&I_{m-r} \\
                                                            \end{array}
                                                          \right)$ and
$A=(0,I_k)$. Clearly, $A\vee B$ contains a $(k+r)$-subspace
$\scriptsize B_1:=\left(\begin{array}{ccc}
 I_r & 0  &0\\
  0& 0 &I_{k} \\
   \end{array}\right)$. Since a basis of $A$  can be extended to a basis of any $(k+r)$-subspace $C\in A\vee B$,
it is easy to prove that every $(k+r)$-subspace containing $A$ and $B$ must be $B_1$. Thus $A\vee B= \{B_1\}$ is a fixed $(k+r)$-subspace.
If ${\rm dim}(A\vee B)=k+r=n$, then ${\rm dim}(A\vee B)=\mathbb{Z}_{p^s}^n$ and hence $A\vee B$ is fixed.

Conversely, suppose $A\vee B$ is a fixed $(k+r)$-subspace. If $k+r=n$, then (ii) holds. Now, we assume that $k+r<n$.
We assert $\scriptsize\rho\left(
\begin{array}{c}
A \\
B \\
\end{array}
\right)={\rm rk}\left(
\begin{array}{c}
A \\
B \\
\end{array}
\right)$. Otherwise, $\scriptsize\rho\left(
\begin{array}{c}
A \\
B \\
\end{array}
\right)>{\rm rk}\left(
\begin{array}{c}
A \\
B \\
\end{array}
\right)$, and hence $0\neq p^{i_r}\in J_{p^s}$ in the matrix $D$. By  elementary row operations of matrix,  it is easy to see that the $(k+r)$-subspace
$\scriptsize B_2:=\left(\begin{array}{ccccc}
 I_{r-1} & 0 &0 & 0 &0\\
 0 & 1 & p^{s-1} & 0 &0\\
  0& 0 & 0& 0& I_k \\
  \end{array}
   \right)$ contains $A$ and $B$. Also, the $(k+r)$-subspace $\small B_1=\left(
                                                            \begin{array}{ccc}
                                                              I_r & 0  &0\\
                                                               0& 0 &I_{k} \\
                                                            \end{array}
                                                          \right)$
  contains $A$ and $B$, a contradiction since $B_1\neq B_2$.
$\qed$

\subsection{Enumeration on subspaces}

\ \ \ \ \
For any non-negative integers  $m, n,q$ with $n\geq m$ and $q\geq 2$, the  {\em Gaussian binomial coefficient} is
$$
{\ n \ \brack m}_q=\prod_{i=1}^m\frac{q^{n+1-i}-1}{q^i-1},
$$
where ${\ n \ \brack 0}_q:=1$. Let ${0\brack m}_q:=0$ (if $m>0$). We have ${n\brack m}_q={n\brack n-m}_q$.

\begin{thm}\label{subspacese3ii5} \ Let  $n> m>k\geq 1$. Then:
\begin{itemize}
\item[{\rm (i)}] The number of $m$-subspaces  of $\mathbb{Z}_{p^s}^n$ is $p^{(s-1)m(n-m)}{n\brack m}_p$.

\item[{\rm (ii)}] In $\mathbb{Z}_{p^s}^n$, the number of $k$-subspaces in a given $m$-subspace is $p^{(s-1)k(m-k)}{m\brack k}_p$.

\item[{\rm (iii)}] In $\mathbb{Z}_{p^s}^n$, the number of $m$-subspaces containing a given $k$-subspace is $p^{(s-1)(m-k)(n-m)}{n-k\brack m-k}_p$.
\end{itemize}
\end{thm}
\proof (i). \ Let $n_{m,s}$ be the number of $m$-subspaces  of $\mathbb{Z}_{p^s}^n$. When $s=1$,
 it is well-known that $n_{m,1}={n\brack m}_p$. Now, we assume $s\geq 2$.
Let $\pi: \mathbb{Z}_{p^s}^{m\times n}\rightarrow \mathbb{Z}_{p}^{m\times n}$ be the natural surjection (\ref{naturalsurjection}).
Suppose $S$ is the set of all $m$-subspaces of $\mathbb{Z}_{p^s}^n$. By Theorem \ref{pipipipi005},
$\pi(S)$ is the set of all $m$-subspaces of $\mathbb{Z}_{p}^n$. Thus $|\pi(S)|={n\brack m}_p$.
Let  $A$ be any fixed $m$-subspace in $\pi(S)$. By Theorem \ref{hgf5tedrdg3}, there is $U_A\in GL_n(\mathbb{Z}_{p})$ such that $A=(0, I_m)U_A$.  It is easy to see that
the preimages $\pi^{-1}(A)=\left\{(X,I_m)U_A: X\in J_{p^s}^{m\times(n-m)}\right\}$.  By $|\pi^{-1}(A)|=p^{(s-1)m(n-m)}$, we obtain $n_{m,s}=p^{(s-1)m(n-m)}{n\brack m}_p$.

(ii). \ By (i), it is clear.

(iii). \ Let $S$ be the set of $m$-subspaces in $\mathbb{Z}_{p^s}^n$ containing a given $k$-subspace $B$ in $\mathbb{Z}_{p^s}^n$.
By Theorem \ref{hgf5tedrdg3}, we can assume that $B=(0, I_k)U$ where $U\in GL_n(\mathbb{Z}_{p^s})$. Then
$$ S=\left\{\left(\begin{array}{cc}
                                 X &0 \\
                                 0& I_k
                               \end{array}
                             \right)U: \mbox{$X$ is any $(m-k)$-subspace of $\mathbb{Z}_{p^s}^{n-k}$}
  \right\}.$$
Using (i), we obtain $|S|=p^{(s-1)(m-k)(n-m)}{n-k\brack m-k}_p$.
$\qed$

\subsection{Dual subspace and arithmetic distance}

\ \ \ \ \
Now, we study dual subspace and arithmetic distance on vector subspace of $\mathbb{Z}_{p^s}^n$.
The dual subspace is an important tool that can simplify some complex calculations.

Let $P$ be an $m$-subspace of $\mathbb{Z}_{p^s}^n$, and let
$$P^\perp=\left\{ y\in\mathbb{Z}_{p^s}^n : \mbox{$y \,^tx=0$ \ for any vector $x\in P$}\right\}.$$
The $P^\perp$ is the set of vectors which are orthogonal to every vector of $P$. Obviously, $P^\perp$ is a linear subset of $\mathbb{Z}_{p^s}^n$.
By Theorem \ref{hgf5tedrdg3}, we can assume that $P=(0, I_m)U$, where $U\in GL_n(\mathbb{Z}_{p^s})$. Thus, it is easy to prove that
 \begin{equation}\label{Dualsubspace02}
\mbox{$ P^\perp=(I_{n-m},0)\,^tU^{-1}$, \  if $P=(0, I_m)U$.}
 \end{equation}
 Hence $P^\perp$ is an $(n-m)$-subspace of $\mathbb{Z}_{p^s}^n$. The subspace $P^\perp$ is called the {\em dual subspace} of $P$.  We have
 \begin{equation}\label{Dualsubspace03}
 {\rm dim}(P)+ {\rm dim}(P^\perp)=n,
 \end{equation}
  \begin{equation}\label{Dualsubspace04}
 (P^\perp)^\perp=P.
 \end{equation}

If $P_1\subseteq P_2$ are two subspaces of $\mathbb{Z}_{p^s}^n$, then
 the definition of $P^\perp$ implies that $P_2^\perp\subseteq P_1^\perp$. It follows from (\ref{Dualsubspace04}) that
\begin{equation}\label{Dualsubspace05}
 \mbox{$P_1\subseteq P_2$ \ $\Longleftrightarrow$ \ $P_2^\perp\subseteq P_1^\perp$.}
 \end{equation}

\begin{lem}\label{Dualsubspace07} \ Let $A$ be an $m$-subspace of $\mathbb{Z}_{p^s}^n$, and let
$\pi: \mathbb{Z}_{p^s}^{m\times n}\rightarrow \mathbb{Z}_{p}^{m\times n}$ be the natural surjection (\ref{naturalsurjection}). Then
\begin{equation}\label{Dualsubspace08}
\pi(A^\perp)=(\pi(A))^\perp.
\end{equation}
\end{lem}
\proof
By Theorem \ref{hgf5tedrdg3}, there is $U\in GL_n(\mathbb{Z}_{p^s})$ such that $A=(0, I_m)U$. From (\ref{Dualsubspace02}) we have
$A^\perp=(I_{n-m},0)\,^tU^{-1}$. Let $U=U_0+U_1p+\cdots +U_{s-1}p^{s-1}$ and  $U^{-1}=Y_0+Y_1p+\cdots +Y_{s-1}p^{s-1}$,
where $U_i, Y_i\in \mathbb{Z}_p^{m\times n}$, $i=0,\ldots, s-1$.  It is clear that $Y_0=U_0^{-1}$.
Then  $^tU^{-1}=\,^tU_0^{-1}+\,^tY_1p+\cdots +\,^tY_{s-1}p^{s-1}$.
Using (\ref{pipipipi003}), we obtain that  $\pi(A)=(0, I_m)U_0$ and $\pi(A^\perp)=(I_{n-m},0)\,^tU_0^{-1}$.
Thus (\ref{Dualsubspace08}) holds.
$\qed$

For two subspaces $A, B$ of $\mathbb{Z}_{p^s}^n$, the {\em arithmetic distance} between $A$ and $B$, denoted by ${\rm ad}(A,B)$,
is defined by
\begin{equation}\label{arithmeticdistance01}
{\rm ad}(A,B)=\rho\left(
\begin{array}{c}
A \\
B \\
\end{array}
\right)-{\rm max}\left\{{\rm dim}(A), {\rm dim}(B)\right\}={\rm dim}(A\vee B)-{\rm max}\left\{{\rm dim}(A), {\rm dim}(B)\right\}.
\end{equation}
Clearly, ${\rm ad}(A, B)\geq 0$ and ${\rm ad}(A, B)={\rm ad}(B, A)$. Moreover, ${\rm ad}(A, B)=0$ $\Leftrightarrow$ $A\subseteq B$ or $B\subseteq A$.

\begin{thm}\label{arithmeticdistance02} \ Let $A$ and $B$ be two subspaces of $\mathbb{Z}_{p^s}^n$. Suppose that  $C$ is a subspace of $\mathbb{Z}_{p^s}^n$
with  ${\rm dim}(C)={\rm min}\left\{{\rm dim}(A), {\rm dim}(B)\right\}$. Then
\begin{equation}\label{arithmeticdistance03}
{\rm ad}(A,B)\leq {\rm ad}(A,C)+{\rm ad}(C,B).
\end{equation}
\end{thm}
\proof
Without loss of generality, we assume that ${\rm dim}(A)={\rm max}\left\{{\rm dim}(A), {\rm dim}(B)\right\}$. Thus ${\rm dim}(B)={\rm min}\left\{{\rm dim}(A), {\rm dim}(B)\right\}$.
By the conditions, ${\rm dim}(C)={\rm dim}(B)$.
By (\ref{dimension-formula02}) and (\ref{arithmeticdistance01}), we have that
${\rm ad}(A,B)={\rm dim}(A\vee B)-{\rm dim}(A)={\rm dim}(B)-{\rm dim}(A\cap B)$,
${\rm ad}(A,C)={\rm dim}(A\vee C)-{\rm dim}(A)={\rm dim}(C)-{\rm dim}(A\cap C)$, and ${\rm ad}(C,B)={\rm dim}(C\vee B)-{\rm dim}(B)={\rm dim}(C)-{\rm dim}(C\cap B)$.
Thus, the inequality (\ref{arithmeticdistance03})  is equivalent to the following inequality:
\begin{equation}\label{arithmeticdistance04}
{\rm dim}(C)\geq {\rm dim}(A\cap C)+{\rm dim}(C\cap B)-{\rm dim}(A\cap B).
\end{equation}

Write $d_1={\rm dim}(A\cap C)$ and $d_2={\rm dim}(C\cap B)$.  Let $A_1$ be a $d_1$-subspace in $A\cap C$ and
$B_1$ a $d_2$-subspace in $C\cap B$. Then $A_1\cap B_1\subseteq A\cap B$ and hence ${\rm dim}(A_1\cap B_1)\leq {\rm dim}(A\cap B)$.
Since $A_1, B_1\subseteq C$, ${\rm dim}(C)\geq {\rm dim}(A_1\vee B_1)$.
Using (\ref{dimension-formula02}) we get ${\rm dim}(A_1\vee B_1)=d_1+d_2-{\rm dim}(A_1\cap B_1)$.
Thus,
$${\rm dim}(C)\geq d_1+d_2-{\rm dim}(A_1\cap B_1)\geq d_1+d_2-{\rm dim}(A\cap B).$$
Therefore, we have proved the inequality (\ref{arithmeticdistance04}).
$\qed$

\begin{thm}\label{arithmeticdistance06} \ Let $n>m$ and let $A$ and $B$ be two $m$-subspaces of $\mathbb{Z}_{p^s}^n$.  Then
\begin{equation}\label{arithmeticdistance07}
{\rm ad}(A,B)= {\rm ad}(A^\perp, B^\perp).
\end{equation}
\end{thm}
\proof
Assume that ${\rm ad}(A,B)=r$.
By (\ref{arithmeticdistance01}) and Theorem \ref{hgf5tedrdg3}, without loss of generality, we assume that $A=(0, I_m)$  and  $B=(D, I_m)$,
where $D={\rm diag}(p^{i_1}, \ldots, p^{i_r}, 0_{m-r, n-m-r})$,  $1\leq r\leq {\rm min}\{m, n-m\}$ and $0\leq i_1\leq \ldots \leq i_r\leq {\rm max}\{s-1, 1\}$.

By (\ref{Dualsubspace03}), both $A^\perp$ and $B^\perp$ are $(n-m)$-subspaces of $\mathbb{Z}_{p^s}^n$.
Applying (\ref{Dualsubspace02}), $A^\perp=(I_{n-m},0)$. Let $B^\perp=(X_1,X_2)$ where $X_1\in\mathbb{Z}_{p^s}^{(n-m)\times (n-m)}$ and
$X_2\in\mathbb{Z}_{p^s}^{(n-m)\times m}$. Since $B^\perp\cdot \,^tB=0$, we have $X_1\,^tD+X_2=0$ and hence $\rho(X_2)=\rho(-X_1\,^tD)\leq r$ because (\ref{rank-02}).
Thus
$${\rm ad}(A^\perp, B^\perp)=\rho\left(
\begin{array}{c}
 A^\perp \\
  B^\perp \\
   \end{array}
    \right)-(n-m)
=\rho\left(
\begin{array}{cc}
 I_{n-m}&0\\
  X_1& X_2 \\
   \end{array}
 \right)-(n-m)=\rho(X_2)\leq r={\rm ad}(A,B).$$
 On the other hand,  using (\ref{Dualsubspace04}) we have similarly ${\rm ad}(A,B)={\rm ad}\left((A^\perp)^\perp, (A^\perp)^\perp)\right)\leq {\rm ad}(A^\perp, B^\perp)$.
Therefore (\ref{arithmeticdistance07}) holds.
$\qed$

\section{Grassmann graphs over $\mathbb{Z}_{p^s}$}
 \setcounter{equation}{0}

\ \ \ \ \ \ In this section, we discuss the basic properties of Grassmann graph  over  $\mathbb{Z}_{p^s}$.
 We will determine the valency of every vertex, the clique number and maximum cliques of the Grassmann graph.

\subsection{Valency, clique number and independence number}

\ \ \ \ \ For any two vertices $A, B$ of the Grassmann graph $G_{p^s}(n,m)$, by Theorem \ref{dimension-formula01} and (\ref{arithmeticdistance01}), we have
\begin{equation}\label{Grassm76hfg31e98n2}
\mbox{$A\sim B$ \ $\Longleftrightarrow$ \ $\rho\left(
\begin{array}{c}
A \\
B \\
\end{array}
\right)={\rm dim}(A\vee B)=m+1$ \ $\Longleftrightarrow$ \ ${\rm ad}(A,B)=1$.}
\end{equation}

By (\ref{arithmeticdistance01}), for vertices $A, B, C$ of $G_{p^s}(n,m)$, $\scriptsize {\rm ad}(A, B)=\rho\left(
\begin{array}{c}
A \\
B \\
\end{array}
\right)-m$  is the {\em arithmetic distance} between $A$ and $B$ in $G_{p^s}(n,m)$. Clearly, ${\rm ad}(A, B)\geq 0$, ${\rm ad}(A, B)=0 \Leftrightarrow A=B$; ${\rm ad}(A, B)={\rm ad}(B, A)$;
${\rm ad}(A, B)\leq {\rm ad}(A, C)+{\rm ad}(C, B)$ by Theorem \ref{arithmeticdistance02}.

\begin{thm}\label{Grassmann000} \ If $A$ and $B$ are two vertices of $G_{p^s}(n,m)$, then
\begin{equation}\label{distance000ewwfs}
d(A,B)={\rm ad}(A, B).
\end{equation}
\end{thm}
\proof
Without loss of generality, we assume that $s\geq 2$.
Let $A,B\in V(G_{p^s}(n,m))$ and ${\rm ad}(A, B)=r$, where  $1\leq r\leq {\rm min}\{m,n-m\}$. By Theorem \ref{hgf5tedrdg3},
we may assume with no loss of generality that  $A=(0,  I_m)$ and $B=(D,  I_m)$,
where $D={\rm diag}(p^{i_1}, \ldots, p^{i_r}, 0_{m-r, n-m-r})$, and $0\leq i_1\leq \ldots \leq i_r\leq s-1$.
Let  $A_j=(D_j,  I_m)$ where $D_j={\rm diag}(p^{i_1}, \ldots, p^{i_j}, 0_{m-j, n-m-j})$, $j=1,\ldots, r-1$.
Then $A\sim A_1\sim \cdots \sim A_{r-1}\sim B$.
Hence  $d(A,B)\leq {\rm ad}(A, B)$. Applying the triangle inequality, it is easy to prove that  $d(A,B)\geq {\rm ad}(A, B)$.
Thus (\ref{distance000ewwfs}) holds.
 $\qed$

By Theorem \ref{subspacese3ii5}(i),  we have
\begin{equation}\label{Grassmann001}
\left|V(G_{p^s}(n,m))\right|=p^{(s-1)m(n-m)}{n\brack m}_p.
\end{equation}

By Theorem \ref{arithmeticdistance06}, (\ref{Grassm76hfg31e98n2}) and (\ref{Dualsubspace02}),  it is easy to see that
$G_{p^s}(n,m)$ is isomorphic to $G_{p^s}(n,n-m)$.
When $m=1$,  $G_{p^s}(n,1)$ a complete graph. Thus,  we always assume that $n\geq 2m \geq 4$  in our discussion on
$G_{p^s}(n,m)$ unless specified otherwise.

\begin{thm}\label{Grassmann004} \ The  $G_{p^s}(n,m)$ is a connected vertex-transitive graph of the valency
\begin{equation}\label{Grassmann005}
r=p^{(s-1)(m-1)}{m\brack 1}_p \left( p^{(s-1)(n-m)}{n-m\brack 1}_p+p^{s(n-m)}-1\right).
\end{equation}
\end{thm}
\proof
Let  $A=(0, I_m)$ be a fixed vertex of of $G_{p^s}(n,m)$.
For any vertex $Y$ of $G_{p^s}(n,m)$,  Theorem \ref{hgf5tedrdg3}  implies that  $Y=(0, I_m)U$ where $U\in GL_n(\mathbb{Z}_{p^s})$. Since the map $X\mapsto XU^{-1}$ is an
 automorphism of $G_{p^s}(n,m)$, $G_{p^s}(n,m)$ is vertex-transitive. Thus $G_{p^s}(n,m)$ is an $r$-regular graph, where $r$ is the same valency of each vertex.
  Without loss of generality, we assume $s\geq 2$.
By (\ref{distance000ewwfs}), it is clear that $G_{p^s}(n,m)$ is connected.

 Let $X$ be any vertex of $G_{p^s}(n,m)$ with $X\sim A$. Write $X=(X_1,X_2)$
where $X_1\in \mathbb{Z}_{p^s}^{m\times (n-m)}$ and  $X_2\in \mathbb{Z}_{p^s}^{m\times m}$.
Since $\scriptsize\rho\left(
\begin{array}{c}
A \\
X \\
\end{array}
\right)=m+1$, $\rho(X_1)=1$. Using Lemma \ref{canonical}, $X$ can be written as
$\scriptsize X=\left(
     \begin{array}{cc}
       \alpha & \beta \\
       0 & X_{22} \\
     \end{array}
   \right)$ where $0\neq\alpha\in\mathbb{Z}_{p^s}^{n-m}$  and $X_{22}\in V(G_{p^s}(m,m-1))$. Thus, (\ref{Grassmann001}) implies that
   $X_{22}$ has $p^{(s-1)(m-1)}{m\brack 1}_p$ different choices. For every fixed choice $X_{22}$, we can assume further that
   $\scriptsize \left(
     \begin{array}{cc}
       \alpha & \beta \\
       0 & X_{22} \\
     \end{array}
   \right)=\left(
     \begin{array}{ccc}
       x & y&0 \\
       0 & 0&I_{m-1} \\
     \end{array}
   \right)$, where $x\in\mathbb{Z}_{p^s}^{n-m}$, $y\in\mathbb{Z}_{p^s}$ and $(x,y)$  is unimodular  with $x\neq 0$, and $A=(0, I_m)$ does not change.

{\bf Case 1}.  $y\in\mathbb{Z}_{p^s}^*$.  We can assume $y=1$. Then for every fixed choice $X_{22}$, $(x, 1)$ with $x\neq 0$ has
 $\left|\mathbb{Z}_{p^s}\right|^{n-m}-1=p^{s(n-m)}-1$ different choices.

{\bf Case 2}. \ $y\in J_{p^s}$.   Then for every fixed choice $X_{22}$, unimodular $(x,y)$ with $x\neq 0$ has
$$\frac{\left|\mathbb{Z}_{p^s}\right|^{n-m}-\left|J_{p^s}\right|^{n-m}}{\left|\mathbb{Z}_{p^s}^*\right|}\left|J_{p^s}\right|=p^{(s-1)(n-m)} {n-m\brack 1}_p$$
 different choices.

 Combination Case 1 with Case 2, for every fixed choice $X_{22}$, unimodular $(x,y)$ with $x\neq 0$ has
 $p^{(s-1)(n-m)}{n-m\brack 1}_p+p^{s(n-m)}-1$ different choices.
 It follows that vertex $X$ with $X\sim A$ has
 $$p^{(s-1)(m-1)}{m\brack 1}_p \left(p^{(s-1)(n-m)}{n-m\brack 1}_p+p^{s(n-m)}-1\right)$$
 different choices. Thus (\ref{Grassmann005}) holds.  $\qed$

For every $(m-1)$-subspace $P$ of $\mathbb{Z}_{p^s}^n$, let $[P\rangle_m$ denote the set of all $m$-subspaces containing $P$, which is called a {\em star}.
For every $(m+1)$-subspace $Q$ of $\mathbb{Z}_{p^s}^n$, let $\langle Q]_m$ denote the set of all $m$-subspaces of $Q$, which is called a {\em top}.
Note that Theorem \ref{dimension-formula01}. It is easy to see that every star or top is a clique of  $G_{p^s}(n,m)$.
In  $G_{p^s}(n,m)$, by Theorem \ref{subspacese3ii5}, we have
\begin{equation}\label{size}
\left|\langle Q]_m \right|=p^{(s-1)m}{m+1\brack 1}_p, \ \ \ \ \ \left|[P\rangle_m \right|=p^{(s-1)(n-m)}{n-m+1\brack 1}_p.
\end{equation}
When $n=2m$,  we have $|[P\rangle_m|=|\langle Q]_m|$. On the other hand, $\left|[P\rangle_m \right|< \left|\langle Q]_m \right|$ if $n<2m$; and
\begin{equation}\label{size00b}
\mbox{$\left|[P\rangle_m \right|> \left|\langle Q]_m \right|$ \ if $n>2m$.}
\end{equation}

Let $P$ and $Q$ be an $m$-subspace and an $(m+1)$-subspace of $\mathbb{Z}_{p^s}^n$, respectively.
By Lemma \ref{hgf5tedrdg3}, there are $U_1, U_2\in GL_n(\mathbb{Z}_{p^s})$ such that $P=(0,I_{m-1})U_1$ and $Q=(0,I_{m+1})U_2$. Thus
\begin{equation}\label{Grassmann008}
[P\rangle_m=\left\{\left(\begin{array}{cc}
                       \alpha & 0 \\
                       0 & I_{m-1} \\
                     \end{array}
                   \right)U_1: \mbox{$\alpha\in\mathbb{Z}_{p^s}^{n-m+1}$ is unimodular}\right\},
\end{equation}
\begin{equation}\label{Grassmann009}
\langle Q]_m=\left\{(0,Y)U_2: Y\in V(G_{p^s}(m+1,m))\right\}.
\end{equation}

When $n=2m$, we define
$$([P\rangle_m)^\perp=\left\{X^\perp: \mbox{$X$ is any $m$-subspace in $[P\rangle_m$}  \right\},$$
$$(\langle Q]_m)^\perp=\left\{X^\perp: \mbox{$X$ is any $m$-subspace in $\langle Q]_m$}  \right\}.$$

 In  $G_{p^s}(2m,m)$, using  (\ref{Dualsubspace05}) and (\ref{size}), it is easy to see that
\begin{equation}\label{Grassmann009bb}
([P\rangle_m)^\perp=\langle P^\perp]_m, \ \ \ (\langle Q]_m)^\perp=[Q^\perp\rangle_m.
\end{equation}

\begin{lem}\label{Grassmann010} \ In  $G_{p^s}(n,m)$,  every star or top is a maximal clique.
\end{lem}
\proof
For every star $[P\rangle_m$, by (\ref{Grassmann008}), we have
\begin{equation}\label{Grassmann011}
[P\rangle_m=\left\{\left(\begin{array}{cc}
                       \alpha & 0 \\
                       0 & I_{m-1} \\
                     \end{array}
                   \right)U_1: \mbox{$\alpha\in\mathbb{Z}_{p^s}^{n-m+1}$ is unimodular}\right\},
\end{equation}
where  $U_1\in GL_n(\mathbb{Z}_{p^s})$ is fixed.
Let $X=(X_1,X_2)U_1\in V(G_{p^s}(n,m))$ and $X\sim Z$ for all $Z\in [P\rangle_m$, where $X_1\in\mathbb{Z}_{p^s}^{m\times (n-m+1)}$
and $X_2\in\mathbb{Z}_{p^s}^{m\times (m-1)}$. Then
$$\rho\left(
        \begin{array}{c}
          Z \\
          X \\
        \end{array}
      \right)=\rho\left(
                    \begin{array}{cc}
                      \alpha & 0 \\
                      0 & I_{m-1} \\
                      X_1 & X_2 \\
                    \end{array}
                  \right)=m-1+\rho\left(
        \begin{array}{c}
          \alpha \\
          X_1 \\
        \end{array}
      \right)=m+1,$$
 for all unimodular $\alpha\in\mathbb{Z}_{p^s}^{n-m+1}$.
Thus $\rho(X_1)=1$. By Lemma \ref{canonical}, $\scriptsize X=\left(\begin{array}{cc}
                       \alpha_1 &  \alpha_2  \\
                        0& X_{2}'\\
                     \end{array}
                   \right)U_1$ where $X'_{2}\in\mathbb{Z}_{p^s}^{(m-1)\times (m-1)}$. Clearly, $X'_2$ is invertible.
It follows that $X\in [P\rangle_m$, and hence $[P\rangle_m$ is a maximal clique.

For every top $\langle Q]_m$, by (\ref{Grassmann009}), we may assume  that
\begin{equation}\label{Grassmann013}
\langle Q]_m=\left\{(0,Y)U_2: Y\in V(G_{p^s}(m+1,m))\right\},
\end{equation}
where  $U_2\in GL_n(\mathbb{Z}_{p^s})$ is fixed.
Let $W\in V(G_{p^s}(n,m))$ and $W\sim Z$ for all $Z\in \langle Q]_m$. By $W\sim (0,I_m)U_2$, we have similarly that
$\scriptsize W=\left(\begin{array}{ccc}
                       \alpha_1 &  \alpha_2 & \alpha_3 \\
                       0 & 0& W_{22}\\
                     \end{array}
                   \right)U_2$, where $\alpha_1\in\mathbb{Z}_{p^s}^{n-m-1}$, $\alpha_2\in\mathbb{Z}_{p^s}$ and $W_{22}\in\mathbb{Z}_{p^s}^{(m-1)\times m}$.
Since  $W_{22}$ has a right inverse, it is easy to see that there is $Y_0\in V(G_{p^s}(m+1,m))$ such that $\scriptsize\rho\left(
        \begin{array}{c}
         (0, W_{22})\\
         Y_0\\
        \end{array}
      \right)=m+1$. By $W\sim (0,Y_0)U_2\in \langle Q]_m$, we get $ \alpha_1=0$. Thus $W\in\langle Q]_m$, and hence $\langle Q]_m$ is a maximal clique.
$\qed$

\begin{thm}\label{Grassmann015} \ {\rm(i)} \ The clique number of $G_{p^s}(n,m)$ (when $n\geq 2m$) is
\begin{equation}\label{Grassmann016}
\omega \left(G_{p^s}(n,m)\right)=p^{(s-1)(n-m)}{n-m+1\brack 1}_p.
\end{equation}

\begin{itemize}
\item[{\rm (ii)}] Let  $\mathcal{M}$ be a maximum clique of $G_{p^s}(n,m)$, and let $\pi: \mathbb{Z}_{p^s}^{m\times n}\rightarrow \mathbb{Z}_{p}^{m\times n}$
be the natural surjection (\ref{naturalsurjection}).
Then when $n>2m$, $\pi(\mathcal{M})$ is a star in $G_{p}(n,m)$. When $n=2m$, $\pi(\mathcal{M})$ is a star or a top in $G_{p}(n,m)$.
\end{itemize}
\end{thm}
\proof
Put $G_{p^s}=G_{p^s}(n,m)$. When $s=1$, this theorem is well-known (cf. \cite{Chow, HuangLvWang, M.Pankov02}). From now on we assume that $s\geq 2$ and $n\geq 2m$.
Write  $\omega=\omega(G_{p^s})$.
Let $\pi: \mathbb{Z}_{p^s}^{m\times n}\rightarrow \mathbb{Z}_{p}^{m\times n}$ be the natural surjection (\ref{naturalsurjection}).

 Let $\mathcal{M}=\{A_1,\ldots, A_{\omega}\}$ be a maximum clique of $G_{p^s}$.
Suppose $\{\pi(A_{i_1}), \ldots, \pi(A_{i_k})\}$ is the set all different elements in $\{\pi(A_1), \ldots, \pi(A_\omega)\}$.
By (\ref{Grassm76hfg31e98n2}) and Theorem \ref{pipipipi005}(i),  $\{\pi(A_{i_1}), \ldots, \pi(A_{i_k})\}$ is a clique of $G_{p}$. Thus
$$\mbox{$k\leq \omega(G_{p})={n-m+1\brack 1}_p$.}$$
Moreover, $\mathcal{M}$ has a  partition into $k$  cliques: $\mathcal{M}=M_1\cup M_2\cup \cdots \cup M_k$, where $M_t$  is a clique with $\pi(M_t)=\pi(A_{i_t})$, $t=1,\ldots,k$.
Thus
$$\omega(G_{p^s})=\sum_{t=1}^k|M_t|.$$

{\rm (i)}. \
Put $n_t=|M_t|$, $t=1,\ldots,k$. Then $M_t=\{\pi(A_{i_t})+B_{t1}, \pi(A_{i_t})+B_{t2}, \ldots, \pi(A_{i_t})+B_{tn_t}\}$, where $B_{tj}\in J_{p^s}^{m\times n}$,
$j=1,\ldots, n_t$.
By Theorem \ref{hgf5tedrdg3}, there is $U_t\in GL_n(\mathbb{Z}_{p})$ such that $\pi(A_{i_t})=(0, I_m)U_t$.
Note that $I_m+B$ is invertible for all $B\in J_{p^s}^{m\times m}$. The matrix representation $\pi(A_{i_t})+B_{tj}$ can be written as
$\pi(A_{i_t})+B_{tj}=(C_{tj}, I_k)U_t$, where $C_{tj}\in J_{p^s}^{m\times(n-m)}$, $j=1,\ldots, n_t$. Therefore,
$$M_t=\left\{(C_{t1}, I_m)U_t, (C_{t2}, I_m)U_t, \ldots, (C_{tn_t}, I_m)U_t \right\}, \ \ t=1,\ldots, k.$$
By (\ref{Grassm76hfg31e98n2}), it is clear that $\left\{C_{t1}, C_{t2}, \ldots, C_{tn_t} \right\}$ is a clique of the bilinear forms graph
$\Gamma(\mathbb{Z}_{p^s}^{m\times(n-m)})$. Since $C_{tj}\in J_{p^s}^{m\times(n-m)}$, from (\ref{bcweu7754332t}) we can write $C_{tj}=D_{tj}p$
where  $D_{tj}\in \mathbb{Z}_{p^{s-1}}^{m\times(n-m)}$,
$j=1,\ldots, n_t$. By Lemma \ref{64gdgd7wrww},  $\left\{D_{t1}, D_{t2}, \ldots, D_{tn_t} \right\}$ is a clique of
$\Gamma(\mathbb{Z}_{p^{s-1}}^{m\times(n-m)})$. By Lemma \ref{gdf5435dt},
$|M_t|=n_t\leq \omega\left(\Gamma(\mathbb{Z}_{p^{s-1}}^{m\times (n-m)})\right)=p^{(s-1)(n-m)}$, $t=1,\ldots,k$. Therefore,
\begin{equation}\label{f432vnvayrr78856}
\omega(G_{p^s})=\sum_{t=1}^k|M_t|\leq kp^{(s-1)(n-m)}\leq  p^{(s-1)(n-m)}\omega(G_{p})=p^{(s-1)(n-m)}{n-m+1\brack 1}_p.
\end{equation}

 By (\ref{size}) and (\ref{f432vnvayrr78856}), every star  is a maximum clique of $G_{p^s}$ and (\ref{Grassmann016}) holds.

{\rm (ii)}. \ Using (\ref{f432vnvayrr78856}) and (\ref{Grassmann016}), it is easy to see that $k={n-m+1\brack 1}_p$, and hence $\pi(\mathcal{M})$ is a maximum clique of $G_{p}$
for every maximum clique $\mathcal{M}$ of $G_{p^s}$. When $n>2m$ (resp. $n=2m$), it is well-known (cf. \cite{HuangLvWang, M.Pankov02}) that every maximum clique of $G_{p}$
is a star (resp. a star or a top). Therefore,   when $n>2m$, $\pi(\mathcal{M})$ is a star of $G_{p}(n,m)$. When $n=2m$, $\pi(\mathcal{M})$ is a star or a top of $G_{p}(n,m)$.
$\qed$

It is a difficult open problem to compute the independence number of Grassmann graph. We can only give the following  estimation of independence number on $G_{p^s}(n,m)$.

\begin{thm}\label{independentset00a} \ Let $n\geq 2m$ and let $G_{p^s}=G_{p^s}(n,m)$. Then
\begin{equation}\label{independentset00b}
 p^{(s-1)(m-1)(n-m)}\alpha (G_{p})\leq\alpha (G_{p^s})\leq p^{(s-1)(m-1)(n-m)}\frac{{n\brack m}_p }{{n-m+1\brack 1}_p}.
\end{equation}
\end{thm}
\proof
Without loss of generality, we assume $s\geq 2$.
Let  $N=\{B_1,\ldots, B_{h}\}$ be a largest independent set of $G_{p}$, where $h=\alpha(G_{p})$.
Applying Theorem \ref{hgf5tedrdg3}, there is $Q_t\in GL_n(\mathbb{Z}_{p})$ such that $B_{t}=(0, I_m)Q_t$, $t=1,\ldots, h$.
Let $\left\{P_{1}, P_{2}, \ldots, P_{r} \right\}$ be a largest independent set of the  bilinear forms graph $\Gamma(\mathbb{Z}_{p^{s-1}}^{m\times(n-m)})$,
where $r=p^{(s-1)(m-1)(n-m)}$ because (\ref{645rerta12bnn}). Put
$$N_t=\left\{(pP_{1}, I_m)Q_t, (pP_{2}, I_m)Q_t, \ldots, (pP_{r}, I_m)Q_t \right\}, \ \ t=1,\ldots, h.$$
Then $|N_t|=r$.
By Lemma \ref{64gdgd7wrww} and (\ref{Grassm76hfg31e98n2}),  $N_t$ is an independent set of $G_{p^s}$, $t=1,\ldots, h$. Note that $(pP_{i}, 0)Q_t$ is a matrix
over $J_{p^s}$, $i=1,\ldots,r$, $t=1,\ldots, h$. For any $X=(pP_{i}, I_m)Q_t\in N_t$ and $Y=(pP_{j}, I_m)Q_u\in N_u$, where $t\neq u$ and $1\leq i,j\leq r$,
by $\scriptsize\rho\left(
\begin{array}{c}
B_t \\
B_u \\
\end{array}
\right)\geq {\rm rk}\left(
\begin{array}{c}
B_t \\
B_u \\
\end{array}
\right)>m+1$ and Lemma \ref{Mc-rank-1},
we have
$$\rho\left(
\begin{array}{c}
X \\
Y \\
\end{array}
\right)=\rho\left(
\begin{array}{c}
(pP_{i}, 0)Q_t+B_t\\
(pP_{j}, 0)Q_u+B_u\\
\end{array}
\right)\geq{\rm rk}\left(
\begin{array}{c}
(pP_{i}, 0)Q_t+B_t\\
(pP_{j}, 0)Q_u+B_u\\
\end{array}
\right)>m+1.
$$
Thus,  $M':=N_1\cup N_2\cup \cdots \cup N_h$ is an independent set of $G_{p^s}$. Consequently
\begin{equation}\label{largestindependent00d}
\alpha(G_{p^s})\geq \sum_{t=1}^h|N_t|=hr= p^{(s-1)(m-1)(n-m)}\alpha(G_{p}).
\end{equation}

If  $G$ is a vertex-transitive graph, then it is well-known (cf. \cite[Lemma 7.2.2]{Godsil}) that
$\alpha(G)\leq \frac{|V(G)|}{\omega(G)}$. Since $G_{p^s}$ is vertex-transitive, (\ref{Grassmann001}) and (\ref{Grassmann016}) imply that
$\alpha (G_{p^s})\leq p^{(s-1)(m-1)(n-m)}\frac{{n\brack m}_p }{{n-m+1\brack 1}_p}$.
$\qed$

\subsection{Maximum cliques}

\begin{thm}\label{Grassmann025} \  When $n>2m$, every  maximum clique of $G_{p^s}(n,m)$ is a star. When $n=2m$, every  maximum clique of $G_{p^s}(n,m)$ is
 either a star or a top.
\end{thm}
\proof
When $s=1$, this theorem is well-known (cf. \cite{Chow,HuangLvWang, M.Pankov02}). From now on we assume that $s\geq 2$ and $n\geq 2m$.
Let $E_{ij}=E_{ij}^{m\times (n-m)}$. Recall that $\Gamma(\mathbb{Z}_{p^s}^{m\times (n-m)})$
 is the  bilinear forms graph on $\mathbb{Z}_{p^s}^{m\times (n-m)}$.

Let $\mathcal{M}$ be a maximum clique of  $G_{p^s}(n,m)$. By Theorem \ref{hgf5tedrdg3}, we may assume with no loss of generality that $A=(0,I_m)\in\mathcal{M}$.
We show that there is a vertex $B\in\mathcal{M}$ such that
$\scriptsize \rho\left(
\begin{array}{c}
A \\
B \\
\end{array}
\right)={\rm rk}\left(
\begin{array}{c}
A \\
B \\
\end{array}
\right)=m+1$. Otherwise, for any vertex $B\in\mathcal{M}$ with $B\sim A$, we have $\scriptsize m+1=\rho\left(
\begin{array}{c}
A \\
B \\
\end{array}
\right)>{\rm rk}\left(
\begin{array}{c}
A \\
B \\
\end{array}
\right)=m$.
It is easy to see that $B=(B_1, I_m)$ where $B_1\in J_{p^s}^{m\times (n-m)}$ with $\rho(B_1)=1$.
  Thus, there exists a clique $\mathcal{N}$ of $\Gamma(\mathbb{Z}_{p^s}^{m\times (n-m)})$ containing  $0$, such that
 $\mathcal{N}\subseteq J_{p^s}^{m\times (n-m)}$ and $\mathcal{M}=\left\{(X, I_m): X\in\mathcal{N}\right\}$.
 Clearly, for every  maximal clique $\mathcal{C}$ containing  $0$ in $\Gamma(\mathbb{Z}_{p^s}^{m\times (n-m)})$, the set $\left\{(X, I_m): X\in\mathcal{C}\right\}$ is a clique of  $G_{p^s}(n,m)$.
Since $\mathcal{M}=\left\{(X, I_m): X\in\mathcal{N}\right\}$ is a maximum clique,  $\mathcal{N}$ must be a maximal clique of $\Gamma(\mathbb{Z}_{p^s}^{m\times (n-m)})$,
but  $\mathcal{N}\subseteq J_{p^s}^{m\times (n-m)}$,  a contradiction to Lemma \ref{gdf5435dt000}.
Therefore, there is vertex $B\in\mathcal{M}$ such that
$\scriptsize \rho\left(
\begin{array}{c}
A \\
B \\
\end{array}
\right)={\rm rk}\left(
\begin{array}{c}
A \\
B \\
\end{array}
\right)=m+1$.  By Theorem \ref{hgf5tedrdg3}, we may assume with no loss of generality that
$$A=(0,I_m), \ \   B=(E_{11}, I_m).$$

Let  $\pi: \mathbb{Z}_{p^s}^{m\times n}\rightarrow \mathbb{Z}_{p}^{m\times n}$ be the natural surjection (\ref{naturalsurjection}).
  For any vertex $X$ in $\pi(\mathcal{M})$, let $\pi^{-1}(X)$ denote the  preimages of $X$ in $\mathcal{M}$, i.e.,
$\pi^{-1}(X)=\left\{Y\in \mathcal{M}: \pi(Y)=X \right\}$.
By Theorem \ref{Grassmann015}(ii), when $n>2m$, $\pi(\mathcal{M})$ is a star in $G_{p}(n,m)$. When $n=2m$, $\pi(\mathcal{M})$ is a star or a top in $G_{p}(n,m)$.

{\bf Case 1}. \  $\pi(\mathcal{M})$ is  a star in $G_{p}(n,m)$.

We will prove that $\mathcal{M}$ is a star in $G_{p^s}(n,m)$.
Since $\pi(A)=A$, $\pi(B)=B$ and $A\cap B=(0, I_{m-1})$,
(\ref{Grassmann008}) implies that
\begin{equation}\label{Grassmann019}
\pi(\mathcal{M})=\left\{\left(\begin{array}{cc}
                       x & 0 \\
                       0 & I_{m-1} \\
                     \end{array}
                   \right): \mbox{$0\neq x\in\mathbb{Z}_{p}^{n-m+1}$}\right\}.
\end{equation}

Let $e_i$ be the $i$-th row of $I_{n-m+1}$, and let
$$A_i=\left(\begin{array}{cc}
                       e_i & 0 \\
                       0 & I_{m-1} \\
                     \end{array}
                   \right)\in \pi(\mathcal{M}), \ \ i=1,\ldots, n-m+1.$$
By (\ref{Grassmann016}),  we have $\left|\pi^{-1}(A_i)\right|=p^{(s-1)(n-m)}$. From Lemma \ref{gdf5435dt}, there is a maximum clique  $\mathcal{C}_i$
of $\Gamma(\mathbb{Z}_{p^{s-1}}^{m\times (n-m)})$ such that
$$
\pi^{-1}(A_i)=\left\{\bordermatrix{
&_{i-1} & _1 & _{ n-m-i+1}& _{m-1} \cr
& \beta_{i1}p &1&\beta_{i2}p &0\cr
& \beta_{i3}p &0& \beta_{i4}p & I_{m-1}  \cr }
: \left(\begin{array}{cc}
                       \beta_{i1}&\beta_{i2}  \\
                       \beta_{i3} & \beta_{i4}\\
                     \end{array}
                   \right)\in \mathcal{C}_i\right\},
$$
$i=1,\ldots, n-m+1$. Here $\scriptsize\left(\begin{array}{c}
                       \beta_{i1}  \\
                       \beta_{i3} \\
                     \end{array}
                   \right)$ or $\scriptsize\left(\begin{array}{c}
                       \beta_{i2}  \\
                       \beta_{i4} \\
                     \end{array}
                   \right)$ may be absent if $i=1$ or $i=n-m+1$.

Let $X_{n-m}\in \pi^{-1}(A_{n-m})$. Then $X_{n-m}\sim A$ and hence $\scriptsize \rho\left(
\begin{array}{c}
A \\
X_{n-m} \\
\end{array}
\right)=m+1$. It follows that $\beta_{n-m,3}=0$. Since $X_{n-m}\sim B$, we have similarly $\beta_{n-m,4}=0$. Therefore, we obtain
\begin{equation}\label{Grassmann022}
\pi^{-1}(A_{n-m})=\left\{\bordermatrix{
&_{n-m-1} & _1 & _{1}& _{m-1} \cr
&\alpha_{n-m,1} &1&\alpha_{n-m,2} &0\cr
&0&0& 0 & I_{m-1}  \cr }
: \alpha_{n-m,1}\in J_{p^s}^{1\times (n-m-1)}, \alpha_{n-m,2}\in J_{p^s}\right\}.
\end{equation}
We can write
$$
\pi^{-1}(A)=\left\{\bordermatrix{
&_{n-m-1} & _{1} & _{1}& _{m-1} \cr
& A_{11}p &  A_{12}p &1 &0\cr
& A_{21}p & A_{22}p & 0 & I_{m-1}  \cr }
: \left(\begin{array}{cc}
                       A_{11}& A_{12}  \\
                       A_{21} & A_{22}\\
                     \end{array}
                   \right)\in \mathcal{C}_{n-m+1}\right\}.
$$
By (\ref{Grassmann022}), we get
$$\left(
                                    \begin{array}{cccc}
                                       A_{11}p &  A_{12}p &1 &0\\
                                      A_{21}p & A_{22}p & 0 & I_{m-1}\\
                                    \end{array}
                                  \right)\sim \left(
                                    \begin{array}{cccc}
                                      0 &1& 0 &0\\
                                      0&0& 0 & I_{m-1}\\
                                    \end{array}
                                  \right),$$
 and hence $\small\rho\left(
                                          \begin{array}{cc}
                                            0 & 1 \\
                                             A_{21}p &  A_{22}p \\
                                          \end{array}
                                        \right)=1$.
Consequently $A_{21}p=0$. Since $\small\left(
                                    \begin{array}{cccc}
                                       A_{11}p &  A_{12}p &1 &0\\
                                      0 & A_{22}p & 0 & I_{m-1}\\
                                    \end{array}
                                  \right)\sim B$, it is easy to see that $A_{22}p=0$. Thus
\begin{equation}\label{Grassmann023}
\pi^{-1}(A)=\left\{\bordermatrix{
&_{n-m} & _{1}& _{m-1} \cr
& \beta &1 &0\cr
& 0& 0 & I_{m-1}  \cr }
: \beta\in J_{p^s}^{1\times (n-m)}\right\}.
\end{equation}

 Let
$$\mbox{$B_\alpha=\left(
\begin{array}{ccc}
 \alpha & 1 & 0 \\
  0 & 0 & I_{m-1} \\
   \end{array}
  \right)\in \pi(\mathcal{M})$, \ \ where $0\neq \alpha\in \mathbb{Z}_p^{n-m}$.}$$
By $\left|\pi^{-1}( B_\alpha)\right|=p^{(s-1)(n-m)}$ and Lemma \ref{gdf5435dt}, there is a maximum clique  $\mathcal{C}_\alpha$
of  $\Gamma(\mathbb{Z}_{p^{s-1}}^{m\times (n-m)})$,  such that
$$\pi^{-1}(B_\alpha)=\left\{\bordermatrix{
&_{n-m} & _{1} & _{m-1} \cr
& \alpha+B_{1}p  &1 &0\cr
& B_{2}p & 0 & I_{m-1}  \cr }
: \left(\begin{array}{c}
                       B_{1} \\
                       B_{2} \\
                     \end{array}
                   \right)\in \mathcal{C}_\alpha\right\}.$$
Let  $e_i'$ be the $i$-th row of $I_{n-m}$.  Since
$$\left(
                                    \begin{array}{ccc}
                                      \alpha+B_{1}p  &1 &0\\
                                      B_{2}p & 0 & I_{m-1}\\
                                    \end{array}
                                  \right)\sim \left(
                                    \begin{array}{ccc}
                                      e_{n-m}' &0 &0\\
                                      0     & 0 & I_{m-1}\\
                                    \end{array}
                                  \right)=A_{n-m}$$
 and (\ref{Grassm76hfg31e98n2}), we  get
$\small \rho\left(
\begin{array}{c}
e_{n-m}'\\
B_{2}p\\
\end{array}
\right)=1$, which implies that $B_{2}p=(0, \beta_2)$ where $\beta_2\in J_{p^s}^{(m-1)\times 1}$.
By (\ref{Grassmann023}), we have
$$\mbox{$\left(
                                    \begin{array}{ccc}
                                      \alpha+B_{1}p  &1 &0\\
                                      B_{2}p & 0 & I_{m-1}\\
                                    \end{array}
                                  \right)\sim \left(
                                    \begin{array}{ccc}
                                      \beta &1 &0\\
                                      0 & 0 & I_{m-1}\\
                                    \end{array}
                                  \right)$ \ \  for all $\beta\in J_{p^s}^{1\times (n-m)}$.}$$
 It follows from (\ref{Grassm76hfg31e98n2}) that
$\small \rho\left(
\begin{array}{c}
\alpha+B_1p-\beta \\
B_{2}p\\
\end{array}
\right)=1$ for all $\beta\in J_{p^s}^{1\times (n-m)}$. Taking $\beta=B_1p$, we get
$\small \rho\left(
\begin{array}{c}
\alpha \\
B_{2}p\\
\end{array}
\right)= \rho\left(
\begin{array}{c}
\alpha \\
(0, \beta_2)\\
\end{array}
\right)=1$. Thus, when $\alpha\neq e_{n-m}'$, we have that $B_{2}p=(0, \beta_2)=0$.
 When $\alpha= e_{n-m}'$,
 $$\left(
                                    \begin{array}{ccc}
                                    e_{n-m}'+B_{1}p  &1 &0\\
                                      B_{2}p & 0 & I_{m-1}\\
                                    \end{array}
                                  \right)\sim \left(
                                    \begin{array}{ccc}
                                      e_{1}' &1 &0\\
                                      0     & 0 & I_{m-1}\\
                                    \end{array}
                                  \right)=B.$$
  By (\ref{Grassm76hfg31e98n2}) we can obtain
$$\rho\left(
\begin{array}{c}
e_{n-m}'+B_{1}p-e_1'\\
B_{2}p\\
\end{array}
\right)=\rho\left(
\begin{array}{c}
e_1'-e_{n-m}'-B_{1}p\\
(0, \beta_2)\\
\end{array}
\right)=1,$$
 which implies that $B_{2}p=(0, \beta_2)=0$. Therefore, we always have $B_{2}p=0$. Then we have proved  that
\begin{equation}\label{Grassmann025c}
\pi^{-1}(B_\alpha)=\left\{\bordermatrix{
&_{n-m} & _{1}& _{m-1} \cr
& \alpha+B_{1}p &1 &0\cr
& 0& 0 & I_{m-1}  \cr }
: B_{1}p\in J_{p^s}^{n-m}\right\}.
\end{equation}

Put
$$\mbox{$D_\alpha=\left(
\begin{array}{ccc}
\alpha & 0 & 0 \\
 0 & 0 & I_{m-1} \\
 \end{array}
 \right)\in \pi(\mathcal{M})$, \ \ where $0\neq \alpha\in\mathbb{Z}_p^{n-m}$.}$$
Then
$$\pi^{-1}(D_\alpha)\subseteq\left\{\bordermatrix{
&_{n-m} & _{1} & _{m-1} \cr
& \alpha+D_{1}p  & D_2p &0\cr
& D_3p & D_4p & I_{m-1}  \cr }
: \left(\begin{array}{cc}
                       D_{1} &D_2 \\
                       D_3&D_4 \\
                     \end{array}
                   \right)\in\mathbb{Z}_{p^{s-1}}^{m\times (n-m+1)}\right\}.$$
Let $\scriptsize\left(
                                    \begin{array}{ccc}
                                      \alpha+D_{1}p  & D_2p &0\\
                                     D_3p & D_4p & I_{m-1}\\
                                    \end{array}
                                  \right)\in \pi^{-1}(D_\alpha)$.
Then  $\scriptsize\left(
                                    \begin{array}{ccc}
                                      \alpha+D_{1}p  & D_2p &0\\
                                     D_3p & D_4p & I_{m-1}\\
                                    \end{array}
                                  \right)\sim A$, and hence
$\scriptsize\rho\left(
\begin{array}{c}
\alpha+D_{1}p \\
D_3p \\
\end{array}
\right)=1$. Since $\alpha+D_{1}p$ is unimodular,  $\scriptsize\left(
\begin{array}{c}
\alpha+D_{1}p \\
D_3p \\
\end{array}
\right)$ has a minimal factorization $\scriptsize\left(
\begin{array}{c}
\alpha+D_{1}p \\
D_3p \\
\end{array}
\right)=\left(
          \begin{array}{c}
            1 \\
          \gamma \\
          \end{array}
        \right)(\alpha+D_{1}p)$, and hence $D_3p=\gamma(\alpha+D_{1}p)$. Therefore, applying elementary row operations of matrix,
we can assume
$$\pi^{-1}(D_\alpha)\subseteq\left\{\bordermatrix{
&_{n-m} & _{1} & _{m-1} \cr
& \alpha+D_{1}p  & D_2p &0\cr
& 0 & D_4'p & I_{m-1}  \cr }
: \left(\begin{array}{cc}
                       D_{1} &D_2 \\
                       0&D_4' \\
                     \end{array}
                   \right)\in \mathbb{Z}_{p^{s-1}}^{m\times (n-m+1)}\right\}.$$

Put $\scriptsize\left(
                                    \begin{array}{ccc}
                                      \alpha+D_{1}p  & D_2p &0\\
                                     0 & D_4'p & I_{m-1}\\
                                    \end{array}
                                  \right)\in \pi^{-1}(D_\alpha)$.
By (\ref{Grassmann025c}), $\scriptsize\left(
                                    \begin{array}{ccc}
                                      \alpha+D_{1}p  & D_2p &0\\
                                     0 & D_4'p & I_{m-1}\\
                                    \end{array}
                                  \right)\sim \left(
                                    \begin{array}{ccc}
                                    \alpha'+B_{1}p &1 &0\\
                                      0 & 0 & I_{m-1}\\
                                    \end{array}
                                  \right)$ for all $B_{1}p\in J_{p^s}^{n-m}$ and $0\neq \alpha'\in \mathbb{Z}_p^{n-m}$. Thus, by (\ref{Grassm76hfg31e98n2}) we get
$\scriptsize\rho\left(
\begin{array}{cc}
 \alpha'+B_{1}p &1 \\
 \alpha+D_{1}p  & D_2p \\
 0 & D_4'p \\
 \end{array}\right)=2$,  for all $0\neq \alpha'\in \mathbb{Z}_p^{n-m}$.
  Clearly, there exists  $0\neq \alpha'\in \mathbb{Z}_p^{n-m}$ such that
 $\scriptsize {\rm rk}\left(\begin{array}{c}
 \alpha'+B_{1}p \\
 \alpha+D_{1}p  \\
 \end{array}\right)=\rho\left(\begin{array}{c}
 \alpha'+B_{1}p  \\
 \alpha+D_{1}p  \\
 \end{array}\right)=2$, which implies that $D_4'p=0$. Therefore,
\begin{equation}\label{Grassmann029}
\pi^{-1}(D_\alpha)=\left\{\bordermatrix{
&_{n-m} & _{1} & _{m-1} \cr
& \alpha+D_{1}p  & D_2p &0\cr
& 0 & 0 & I_{m-1}  \cr }
: (D_1,D_2)\in\mathbb{Z}_{p^{s-1}}^{n-m+1}\right\}.
\end{equation}

Let $P=A\cap B=(0,0,I_{m-1})$. By (\ref{Grassmann023})-(\ref{Grassmann029}), we have that $\mathcal{M}=\pi^{-1}(\pi(\mathcal{M}))=[P\rangle_m$ is a star  in $G_{p^s}(n,m)$.

{\bf Case 2}. \  $\pi(\mathcal{M})$ is  a top in $G_{p}(n,m)$ with $n=2m$.

We show that $\mathcal{M}$ is a top in $G_{p^s}(n,m)$ as follows. Define
$\mathcal{M}^\perp=\left\{X^\perp: \mbox{$X$ is any $m$-subspace in $\mathcal{M}$} \right\}$.
By Lemma \ref{Dualsubspace07}, we have $\pi(\mathcal{M}^\perp)=(\pi(\mathcal{M}))^\perp$.
 It follows from (\ref{Grassmann009bb}) that $(\pi(\mathcal{M}))^\perp$ is  a star in $G_{p}(n,m)$. By the Case 1, $\mathcal{M}^\perp$ is a star in $G_{p^s}(n,m)$.
Applying (\ref{Dualsubspace04}) and (\ref{Grassmann009bb}), $\mathcal{M}=(\mathcal{M}^\perp)^\perp$ is a top in $G_{p^s}(n,m)$.

Combination Case 1 with Case 2, when $n>2m$, every  maximum clique of $G_{p^s}(n,m)$ is a star. When $n=2m$, every  maximum clique of $G_{p^s}(2m,m)$ is
 either a star or a top.
$\qed$

\section{Automorphisms of Grassmann graph over $\mathbb{Z}_{p^s}$}
 \setcounter{equation}{0}

\ \ \ \ \ Recall that an {\em endomorphism} of a graph $G$ is an adjacency preserving map from $V(G)$ to itself.
An {\em automorphism} of a graph $G$ is an adjacency preserving bijective map from $V(G)$ to itself.
In the algebraic graph theory \cite{Godsil}, the characterization of graph automorphism  is an important problem.
In this section, we determine the automorphisms of  $G_{p^s}(n,m)$, and our main result is as follows.

\begin{thm}\label{automorphism014} \ Let $n\geq 2m\geq 4$, and let $\varphi$ be an automorphism of $G_{p^s}(n,m)$. Then either there exists
 $U\in GL_n(\mathbb{Z}_{p^s})$ such that
\begin{equation}\label{automorphism015}
\varphi(X)=XU, \ \ X\in V(G_{p^s}(n,m));
\end{equation}
or $n=2m$ and there exists  $U\in GL_n(\mathbb{Z}_{p^s})$ such that
\begin{equation}\label{automorphism016}
\varphi(X)=(XU)^\perp, \ \ X\in V(G_{p^s}(2m,m)).
\end{equation}
\end{thm}

Theorem \ref{automorphism014} also has important significance for the {\em geometry of matrices} \cite{Huang-book,Wbook}, and it is also the fundamental theorem of the projective geometry
of rectangular matrices over $\mathbb{Z}_{p^s}$. To prove Theorem \ref{automorphism014}, we need the following knowledge and lemmas.

In $G_{p^s}(n,m)$ (where $n>m\geq 1$), two vertices $A,B$ are said to be {\em McCoy adjacent} ({\em Mc-adjacent} for short),
denoted by $A\stackrel{mc}{\sim}B$, if  $\scriptsize \rho\left(
\begin{array}{c}
A \\
B \\
\end{array}
\right)={\rm rk}\left(
\begin{array}{c}
A \\
B \\
\end{array}
\right)=m+1$.  Two Mc-adjacent vertices are adjacent but  not vice versa.
An endomorphism  $\varphi$ of $G_{p^s}(n,m)$ is called {\em to preserve Mc-adjacency} if  $A\stackrel{mc}{\sim}B$ implies that $\varphi(A)\stackrel{mc}{\sim} \varphi(B)$
for any two vertices $A$ and $B$.

\begin{lem}\label{automorphism000} \  Let $n>m\geq 1$, and let $A,B$ be two distinct vertices of $G_{p^s}(n,m)$ with $d(A,B)=k$.
Then there are vertices $A_1,\ldots, A_{2k-1}$ of $G_{p^s}(n,m)$ such that
$A\stackrel{mc}{\sim}A_1 \stackrel{mc}{\sim} \cdots \stackrel{mc}{\sim} A_{2k-1} \stackrel{mc}{\sim} B$.
\end{lem}
\proof
Let $E_{ij}=E_{ij}^{m\times (n-m)}$.
Suppose that $A\sim B$. By Theorem \ref{hgf5tedrdg3}, there is $U\in GL_n(\mathbb{Z}_{p^s})$ such that
$A=(0, I_m)U$  and  $B=(p^iE_{11}, I_m)U$, where  $0\leq i\leq {\rm max}\{s-1,1\}$. Put $C=(p^iE_{11}+E_{12}, I_m)U$. Then
$A\stackrel{mc}{\sim}C \stackrel{mc}{\sim} B$.
Now, let $d(A,B)=k>1$. Then there are vertices $B_1,\ldots, B_{k-1}$  such that $A\sim B_1\sim \cdots \sim B_{k-1}\sim B$.
Thus, there are vertices $C_1,\ldots, C_{k}$  such that
$A\stackrel{mc}{\sim}C_1 \stackrel{mc}{\sim}B_1 \cdots \stackrel{mc}{\sim} B_{k-1} \stackrel{mc}{\sim} C_{k}\stackrel{mc}{\sim} B$.
$\qed$

\begin{lem}\label{automorphism001} \   Let $n>m\geq 2$. Then $A$ is a vertex of $G_{p^s}(n,m)$ if and only if  there are  two vertices $P_1,P_2\in V(G_{p^s}(n,m-1))$
 such that $P_1\stackrel{mc}{\sim}P_2$ and $A=\{A\}=P_1\vee P_2$.
\end{lem}
\proof
Let $A$ be a vertex of $G_{p^s}(n,m)$. By Theorem \ref{hgf5tedrdg3}, we may assume with no loss of generality that $A=(0, I_m)$.
Put $P_1=(0,I_{m-1}), P_2=(0,I_{m-1}, 0_{m,1})\in V(G_{p^s}(n,m-1))$. Then $\scriptsize \rho\left(
\begin{array}{c}
P_1 \\
P_2 \\
\end{array}
\right)={\rm rk}\left(
\begin{array}{c}
P_1 \\
P_2 \\
\end{array}
\right)=m$, and hence $P_1\stackrel{mc}{\sim}P_2$. By Theorem \ref{08ujivxwfal879u}(ii), we get $A=\{A\}=P_1\vee P_2$.

Conversely, suppose there are vertices $P_1,P_2\in V(G_{p^s}(n,m-1))$ such that $P_1\stackrel{mc}{\sim}P_2$ and $A=\{A\}=P_1\vee P_2$.
Then $A$ is a fixed $m$-subspace and hence $A$ is a vertex of $G_{p^s}(n,m)$.
$\qed$

\begin{lem}\label{automorphism002} \  Let $[P_1\rangle_m, [P_2\rangle_m$ be two distinct stars in $G_{p^s}(n,m)$ where $n\geq 2m\geq 4$. Then:
\begin{itemize}
\item[{\rm (i)}]   $[P_1\rangle_m\cap[P_2\rangle_m\neq \emptyset$ if and only if $P_1\sim P_2$ (in $G_{p^s}(n,m-1)$). Moreover,
if $[P_1\rangle_m\cap[P_2\rangle_m\neq \emptyset$,
then  $[P_1\rangle_m\cap[P_2\rangle_m =P_1\vee P_2$.

\item[{\rm (ii)}]  In $G_{p^s}(n,m)$, $\left| [P_1\rangle_m\cap[P_2\rangle_m \right|=1$ if and only if $P_1\stackrel{mc}{\sim}P_2$ (in $G_{p^s}(n,m-1)$).

\end{itemize}
\end{lem}
\proof
(i). \ Assume that $[P_1\rangle_m\cap[P_2\rangle_m\neq \emptyset$. Let vertex $C\in [P_1\rangle_m\cap[P_2\rangle_m$. Then $P_1,P_2\subset C$.
It follows from  Remark \ref{dimension-formula08} that  ${\rm dim}(P_1\vee P_2)=m$ and $C\in P_1\vee P_2$.
Thus $P_1\sim P_2$ (in $G_{p^s}(n,m-1)$) and  $[P_1\rangle_m\cap[P_2\rangle_m\subseteq P_1\vee P_2$.

Conversely, suppose that $P_1\sim P_2$ (in $G_{p^s}(n,m-1)$). Then (\ref{Grassm76hfg31e98n2}) implies that ${\rm dim}(P_1\vee P_2)=m$.
By Remark \ref{dimension-formula08},  $P_1\vee P_2\subseteq [P_1\rangle_m\cap[P_2\rangle_m$
and hence $[P_1\rangle_m\cap[P_2\rangle_m\neq \emptyset$. Therefore, $[P_1\rangle_m\cap[P_2\rangle_m\neq \emptyset$ if and only if $P_1\sim P_2$ (in $G_{p^s}(n,m-1)$).

Now, we assume $[P_1\rangle_m\cap[P_2\rangle_m\neq \emptyset$. Then $\scriptsize \rho\left(
\begin{array}{c}
P_1 \\
P_2 \\
\end{array}
\right)=m={\rm dim}(P_1\vee P_2)$ and  $[P_1\rangle_m\cap[P_2\rangle_m\subseteq P_1\vee P_2$.
By Remark \ref{dimension-formula08} again,  we have $P_1\vee P_2\subseteq [P_1\rangle_m\cap[P_2\rangle_m$. Thus $[P_1\rangle_m\cap[P_2\rangle_m=P_1\vee P_2$.

(ii). \ Suppose that $\left| [P_1\rangle_m\cap[P_2\rangle_m \right|=1$ in $G_{p^s}(n,m)$. By (i), $P_1\vee P_2$ is a fixed $m$-subspace.
By Theorem \ref{08ujivxwfal879u}(ii), we have
$\scriptsize \rho\left(
\begin{array}{c}
P_1 \\
P_2 \\
\end{array}
\right)={\rm rk}\left(
\begin{array}{c}
P_1 \\
P_2 \\
\end{array}
\right)=m$, and hence $P_1\stackrel{mc}{\sim}P_2$ (in $G_{p^s}(n,m-1)$).

Conversely, suppose that $P_1\stackrel{mc}{\sim}P_2$ (in $G_{p^s}(n,m-1)$). Using (i) and Theorem \ref{08ujivxwfal879u}(ii), we obtain
$\left| [P_1\rangle_m\cap[P_2\rangle_m \right|=|P_1\vee P_2|=1$.
$\qed$

\begin{lem}\label{automorphism001b} \ In $G_{p^s}(n,m)$, if $[P\rangle_m\cap\langle Q]_m\neq \emptyset$, then $\left| [P\rangle_m\cap\langle Q]_m \right|=p^{s-1}(p+1)$.
\end{lem}
\proof
When $s=1$, this lemma is known by \cite[Lemma 2.1]{HuangLvWang}. From now on we assume that $s\geq 2$.
Suppose $[P\rangle_m\cap\langle Q]_m\neq \emptyset$. Then $P\subset Q$. By Lemma \ref{canonical},
without loss of generality, we may assume that $P=(0, I_{m-1})$. Thus $\scriptsize Q=\left(
                                                                           \begin{array}{cc}
                                                                             Q_1 & 0 \\
                                                                             0 & I_{m-1} \\
                                                                           \end{array}
                                                                         \right)$, where  $Q_1\in\mathbb{Z}_{p^s}^{2\times(n-m+1)}$ has a right inverse.
 Using appropriate elementary operations of matrix, we may assume with no loss of generality that $Q=(0, I_{m+1})$.

Let vertex $A\in [P\rangle_m\cap\langle Q]_m$. Then $P\subset A\subset Q$. By (\ref{Grassmann008}) and (\ref{Grassmann009}),  we have
$\scriptsize A=\left(\begin{array}{ccc}
0& \alpha_1 & 0 \\
 0 & 0 & I_{m-1} \\
 \end{array}
 \right)$ where $\alpha_1\in \mathbb{Z}_{p^{s}}^2$ is unimodular. Since the number of unimodular vectors in $\mathbb{Z}_{p^{s}}^2$ is
$\frac{\left|\mathbb{Z}_{p^s}\right|^{2}-\left|J_{p^s}\right|^{2}}{\left|\mathbb{Z}_{p^s}^*\right|}=p^{s-1}(p+1)$,
we obtain $\left| [P\rangle_m\cap\langle Q]_m \right|=p^{s-1}(p+1)$.
$\qed$

\begin{lem}\label{automorphism003} \ Let $n\geq 2m\geq 4$. Suppose  $\varphi$ is an automorphism of $G_{p^s}(n,m)$ and $\varphi$ maps stars to stars.
Let $\varphi([X\rangle_m)=[X'\rangle_m$, $X\in V(G_{p^s}(n,m-1))$.
 Define the map $\varphi_1: V(G_{p^s}(n,m-1))\rightarrow V(G_{p^s}(n,m-1))$ by
$\varphi_1(X)=X'$. Then:
\begin{itemize}
\item[{\rm (i)}]
The map $\varphi_1$ is an automorphism of $G_{p^s}(n,m-1)$,  $\varphi$ is uniquely determined by $\varphi_1$, and
\begin{equation}\label{automorphism005}
\mbox{$\varphi(X\vee Y)=\varphi_1(X)\vee \varphi_1(Y)$,  \ for all $X,Y\in V(G_{p^s}(n,m-1))$ with $X\stackrel{mc}{\sim} Y$.}
\end{equation}

\item[{\rm (ii)}]  Let $X,Y, W\in V(G_{p^s}(n,m-1))$. Then  $X\stackrel{mc}{\sim} Y$ if and only if  $\varphi_1(X)\stackrel{mc}{\sim} \varphi_1(Y)$.
Moreover,  $X\stackrel{mc}{\sim} Y$ with $W\subset X\vee Y$, if and only if $\varphi_1(X)\stackrel{mc}{\sim} \varphi_1(Y)$ with $\varphi_1(W)\subset \varphi_1(X)\vee \varphi_1(Y)$.
\end{itemize}
\end{lem}
\proof
(i). \ By Lemma \ref{automorphism002}(i), it is easy to see that $\varphi_1$ is an injective endomorphism of $G_{p^s}(n,m-1)$. Thus
 Thus  $\varphi_1$ is an automorphism of $G_{p^s}(n,m-1)$.

Put $X,Y\in V(G_{p^s}(n,m-1))$ with $X\stackrel{mc}{\sim} Y$. From Theorem \ref{08ujivxwfal879u}(ii) we know that $X\vee Y$ is a fixed $m$-subspace.
By Lemma \ref{automorphism002}(i), $X\vee Y=[X\rangle_m\cap[Y\rangle_m$. Thus $\varphi(X\vee Y)=[\varphi_1(X)\rangle_m\cap[\varphi_1(Y)\rangle_m$, and hence
$\left| [\varphi_1(X)\rangle_m\cap[\varphi_1(Y)\rangle_m \right|=1$. It follows from Lemma \ref{automorphism002} that $\varphi_1(X)\stackrel{mc}{\sim} \varphi_1(Y)$ and
$$\varphi(X\vee Y)=[\varphi_1(X)\rangle_m\cap[\varphi_1(Y)\rangle_m=\varphi_1(X)\vee \varphi_1(Y).$$
Therefore,  (\ref{automorphism005}) holds and $\varphi_1$ preserves the Mc-adjacency in $G_{p^s}(n,m-1)$.

 Let $A$ be a vertex of $G_{p^s}(n,m)$. By Lemma \ref{automorphism001}, there are two vertices $P_1,P_2\in V(G_{p^s}(n,m-1))$ such that $P_1\stackrel{mc}{\sim}P_2$ and
$A=P_1\vee P_2$. By (\ref{automorphism005}),  $\varphi$  is uniquely determined by $\varphi_1$.

(ii). \ Let $\varphi^{-1}$ and $\varphi_1^{-1}$ be the inverse maps of $\varphi$ and $\varphi_1$, respectively. Then $\varphi^{-1}$ and $\varphi_1^{-1}$
are also two graph automorphisms. Since $\varphi([X\rangle_m)=[\varphi_1(X)\rangle_m$ for all $X\in V(G_{p^s}(n,m-1))$,
we have $\varphi^{-1}([Z\rangle_m)=[\varphi_1^{-1}(Z)\rangle_m$ for all $Z\in V(G_{p^s}(n,m-1))$. By (i), both $\varphi_1$ and  $\varphi_1^{-1}$ preserve the Mc-adjacency in $G_{p^s}(n,m-1)$.
 Therefore, for any $X,Y\in V(G_{p^s}(n,m-1))$, it is easy to see that $X\stackrel{mc}{\sim} Y$ if and only if
$\varphi_1(X)\stackrel{mc}{\sim} \varphi_1(Y)$.

Suppose that $X,Y, W\in V(G_{p^s}(n,m-1))$, $X\stackrel{mc}{\sim} Y$ and $W\subset X\vee Y$. Without losing generality, we assume that $W\neq X$ and $W\neq Y$.
Then   ${\rm dim}(W\vee Y)={\rm dim}(X\vee Y)=m$.
 Hence $W\sim Y$. Similarly, $W\sim X$. By (i), we have that $\varphi_1(X)\stackrel{mc}{\sim} \varphi_1(Y)$, $\varphi_1(W)\sim\varphi_1(Y)$ and $\varphi_1(W)\sim\varphi_1(X)$.
 By Theorem \ref{08ujivxwfal879u}(ii), $X\vee Y$ and $\varphi_1(X)\vee\varphi_1(Y)$ are two
 fixed vertices of $G_{p^s}(n,m)$, and $X\vee Y\in W\vee Y$.

 Using Lemma \ref{automorphism002}(i),
$[X\rangle_m\cap[Y\rangle_m \in [W\rangle_m\cap[Y\rangle_m$, and hence
$$\varphi([X\rangle_m\cap[Y\rangle_m)=\varphi([X\rangle_m)\cap\varphi([Y\rangle_m) \in \varphi([W\rangle_m\cap[Y\rangle_m)=\varphi([W\rangle_m)\cap\varphi([Y\rangle_m).$$
By Lemma \ref{automorphism002}(i) again, we have $\varphi_1(X)\vee \varphi_1(Y)\in \varphi_1(W)\vee \varphi_1(Y)$, it follows from Remark \ref{dimension-formula08} that
 $\varphi_1(W)\subset \varphi_1(X)\vee \varphi_1(Y)$.
Thus, if $X\stackrel{mc}{\sim} Y$ with $W\subset X\vee Y$, then $\varphi_1(X)\stackrel{mc}{\sim} \varphi_1(Y)$ with $\varphi_1(W)\subset \varphi_1(X)\vee \varphi_1(Y)$.
Similarly, considering  $\varphi_1^{-1}$,  we have that  $\varphi_1(X)\stackrel{mc}{\sim} \varphi_1(Y)$ with $\varphi_1(W)\subset \varphi_1(X)\vee \varphi_1(Y)$ implies that
$X\stackrel{mc}{\sim} Y$ with $W\subset X\vee Y$.
$\qed$

\begin{lem}\label{automorphism008} \ Let $\varphi$ be an automorphism of $G_{p^s}(2m,m)$ (where $m\geq 2$). If $\varphi$ maps some star to a star,
then $\varphi$ maps stars to stars. If $\varphi$ maps some star to a top, then $\varphi$ maps stars to tops.
\end{lem}
\proof
{\bf Case 1}. \ $\varphi$ maps some star to a star.

Then there is a star $[P_1\rangle_m$  such that $\varphi([P_1\rangle_m)=[P_1'\rangle_m$.
We prove that $\varphi$ maps stars to stars as follows.

Let $[P_2\rangle_m$ be any star in $G_{p^s}(2m,m)$ with $P_1\stackrel{mc}{\sim}P_2$ (in $G_{p^s}(2m,m-1)$).
Note that $\varphi$ carries  maximum cliques of $G_{p^s}(2m,m)$ to  maximum cliques.
By Theorem \ref{Grassmann025},
$\varphi([P_2\rangle_m)$ is a star or a top. Suppose $\varphi([P_2\rangle_m)$ is a top. Write $\varphi([P_2\rangle_m)=\langle Q]_m$. Then
$$\varphi([P_1\rangle_m\cap[P_2\rangle_m)=\varphi([P_1\rangle_m)\cap\varphi([P_2\rangle_m)=[P_1'\rangle_m\cap \langle Q]_m.$$
By Lemma \ref{automorphism002}(ii), we have $\left|\varphi([P_1\rangle_m\cap[P_2\rangle_m)\right|=\left| [P_1\rangle_m\cap[P_2\rangle_m \right|=1$.
But by Lemma \ref{automorphism001b} we get $\left|[P_1'\rangle_m\cap \langle Q]_m\right|=p^{s-1}(p+1)>1$, a contradiction.
Therefore, $\varphi([P_2\rangle_m)$ is a star for any star $[P_2\rangle_m$ in $G_{p^s}(2m,m)$ with $P_1\stackrel{mc}{\sim}P_2$.

Now, let $[P\rangle_m$ be any star in $G_{p^s}(2m,m)$ with $d(P_1,P)=k>1$. By Lemma \ref{automorphism000},
 there are vertices $A_1,\ldots, A_{2k-1}$ of $G_{p^s}(2m,m-1)$ such that
$P_1\stackrel{mc}{\sim}A_1 \stackrel{mc}{\sim} \cdots \stackrel{mc}{\sim} A_{2k-1} \stackrel{mc}{\sim} P$.
Applying the above result, all $\varphi([A_1\rangle_m), \ldots, \varphi([A_{2k-1}\rangle_m), \varphi([P\rangle_m)$ are stars.
Therefore, $\varphi$ maps stars to stars.

{\bf Case 2}. \ $\varphi$ maps some star to a top.

Let $\psi(X)=(\varphi(X))^\perp$,  $X\in V(G_{p^s}(2m,m))$.
By (\ref{Dualsubspace02}) and (\ref{Dualsubspace04}), it is clear that $\psi$ is a surjection from $V(G_{p^s}(2m,m))$ to itself.  By (\ref{arithmeticdistance07}),
$\psi$ is also an automorphism of $G_{p^s}(2m,m)$. Using (\ref{Grassmann009bb}), $\psi$ maps some star to a star. It follows from  Step 1 that
$\psi$ maps stars to stars. Applying (\ref{Grassmann009bb}) again, $\varphi$ maps stars to tops.
$\qed$

Let $R$ be a commutative local ring and $n\geq 3$. A $1$-subspace of $R^n$ is called a {\em line},
and a $2$-subspace of $R^n$ is called a {\em plane}. Let ${\rm P}(R^n)$ be the set of lines of $R^n$, and ${\rm P}(R^n)$ is called the {\em projective space} of $R^n$.
A bijective map $f: {\rm P}(R^n)\rightarrow {\rm P}(R^n)$ is called a {\em projectivity} or {\em collineation} if $f$ satisfies the following properties:
(a) \ For lines $L_1$ and $L_2$, $L_1+L_2$ (where $L_1+L_2$ is the direct sum of modules) is a plane if and only if $f(L_1)+f(L_2)$ is a plane.
(b) \ Suppose $L_1+L_2$ is a plane and $L$ is a line. Then  $L\subset L_1+L_2$ if and only if $f(L)\subset f(L_1)+f(L_2)$. Thus, a projectivity is a bijection between
lines which preserves planes (cf. \cite{mcdonald2}).

\begin{rem}\label{localringplane} \
 When $R=\mathbb{Z}_{p^s}$, for lines $L_1$ and $L_2$, $L_1+L_2$ is a plane if and only if $L_1\stackrel{mc}{\sim} L_2$ in $V(G_{p^s}(n,1))$.
 Moreover, if $L_1\stackrel{mc}{\sim} L_2$, then $L_1+L_2=L_1\vee L_2$ by Lemma \ref{08ujivxwfal879u}(ii).
  \end{rem}

We have the following fundamental theorem  of projective geometry \cite{mcdonald2} over $R$ :

\begin{thm}\label{automorphism010}{\rm(see \cite[Theorem I.13]{mcdonald2})} \ Let $R$ be a commutative local ring and $n\geq 3$. Suppose that
$f: {\rm P}(R^n)\rightarrow {\rm P}(R^n)$ is  a  projectivity. Then, there is a semi-linear bijection $\phi: R^n\rightarrow R^n$ such that
$f={\rm P}(\phi)$.
\end{thm}

\noindent{\bf Proof of Theorem \ref{automorphism014}.} \
Let $\varphi$ be an automorphism of $G_{p^s}(n,m)$ where $n\geq 2m\geq 4$.
Since $\varphi$ carries  maximum cliques to maximum cliques, Theorem \ref{Grassmann025} and Lemma \ref{automorphism008} imply that
$\varphi$ maps stars to stars, or $\varphi$ maps stars to tops with $n=2m$.

{\bf Case 1}. \  $\varphi$ maps stars to stars.

Let
$$
\mbox{$\varphi([X\rangle_m)=[X'\rangle_m$,  \ \ $X\in V(G_{p^s}(n,m-1))$.}
$$
Define the map $\varphi_1: V(G_{p^s}(n,m-1))\rightarrow V(G_{p^s}(n,m-1))$ by $\varphi_1(X)=X'$. By Lemma \ref{automorphism003},
$\varphi_1$ is an automorphism of $G_{p^s}(n,m-1)$, $\varphi$ is uniquely determined by $\varphi_1$, and
\begin{equation}\label{automorphism019}
\mbox{$\varphi(X\vee Y)=\varphi_1(X)\vee \varphi_1(Y)$,  \ for all $X,Y\in V(G_{p^s}(n,m-1))$ with $X\stackrel{mc}{\sim} Y$.}
\end{equation}
Moreover,  for $X,Y, W\in V(G_{p^s}(n,m-1))$,  $X\stackrel{mc}{\sim} Y$ if and only if  $\varphi_1(X)\stackrel{mc}{\sim} \varphi_1(Y)$;
$X\stackrel{mc}{\sim} Y$ with $W\subset X\vee Y$ if and only if $\varphi_1(X)\stackrel{mc}{\sim} \varphi_1(Y)$ with $\varphi_1(W)\subset \varphi_1(X)\vee \varphi_1(Y)$.

Note that every maximum clique of $G_{p^s}(n,m-k)$ ($1\leq k\leq m-1$) is a star. The $\varphi_1$ maps stars to stars.
Put $\varphi_0=\varphi$. Similarly, we can define an automorphism
$\varphi_k: V(G_{p^s}(n,m-k))\rightarrow V(G_{p^s}(n,m-k))$, such that $\varphi_{k-1}$ is uniquely determined by $\varphi_k$ and
\begin{equation}\label{automorphism021}
\mbox{$\varphi_{k-1}(X\vee Y)=\varphi_k(X)\vee \varphi_k(Y)$,  \ for all $X,Y\in V(G_{p^s}(n,m-k))$ with $X\stackrel{mc}{\sim} Y$,}
\end{equation}
 $k=1,\ldots, m-1$. Moreover, for all $X,Y, W\in V(G_{p^s}(n,m-k))$,  $X\stackrel{mc}{\sim} Y$ if and only if  $\varphi_k(X)\stackrel{mc}{\sim} \varphi_k(Y)$;
$X\stackrel{mc}{\sim} Y$ with $W\subset X\vee Y$ if and only if $\varphi_k(X)\stackrel{mc}{\sim} \varphi_k(Y)$ with $\varphi_k(W)\subset \varphi_k(X)\vee \varphi_k(Y)$.
 $k=1,\ldots, m-1$.

Let $R=\mathbb{Z}_{p^s}$.
Clearly,   ${\rm P}(R^n)=V(G_{p^s}(n,1))$ and $\varphi_{m-1}$ is a bijective map from ${\rm P}(R^n)$ to itself.
By Remark \ref{localringplane}, $\varphi_{m-1}: {\rm P}(R^n)\rightarrow {\rm P}(R^n)$ is  a  projectivity.
By Theorem \ref{automorphism010}, there is a semi-linear bijection $\phi: R^n\rightarrow R^n$ such that $\varphi_{m-1}={\rm P}(\phi)$.
Since $1$  generates $\mathbb{Z}_{p^s}$, every ring automorphism of $\mathbb{Z}_{p^s}$ is the identity mapping.
Consequently, every  semi-linear bijection $\phi: R^n\rightarrow R^n$ is a linear bijection.
Therefore,  there exists  $U\in GL_n(\mathbb{Z}_{p^s})$ such that
\begin{equation}\label{automorphism024}
\mbox{$\varphi_{m-1}(X)=XU$, \  $X\in V(G_{p^s}(n,1))$.}
\end{equation}

For any $X\in V(G_{p^s}(n,2))$, by Lemma \ref{automorphism001}, there are $X_1,X_2\in V(G_{p^s}(n,1))$ such that $\scriptsize X=X_1\vee X_2=\left(\begin{array}{c}
                                                X_1 \\
                                                X_2 \\
                                              \end{array}
                                            \right)$ with $X_1\stackrel{mc}{\sim} X_2$.
Since $\varphi_{m-1}(X_1)\stackrel{mc}{\sim} \varphi_{m-1}(X_2)$, from (\ref{automorphism021}) and (\ref{automorphism024}), we obtain
$$\varphi_{m-2}(X)=\varphi_{m-1}(X_1)\vee\varphi_{m-1}(X_2)=\left(\begin{array}{c}
                                                \varphi(X_1) \\
                                                 \varphi(X_2) \\
                                              \end{array}
                                            \right)=\left(\begin{array}{c}
                                                X_1 \\
                                                 X_2 \\
                                              \end{array}
                                            \right)U.$$
 Thus
\begin{equation}\label{automorphism026}
\mbox{$\varphi_{m-2}(X)=XU$,  \  $X\in V(G_{p^s}(n,2))$.}
\end{equation}

For any $X\in V(G_{p^s}(n,3))$, by Lemma \ref{automorphism001}, there are $X_1,X_2\in V(G_{p^s}(n,2))$ such that $\scriptsize X=X_1\vee X_2$ with
$X_1\stackrel{mc}{\sim} X_2$. By Theorems \ref{hgf5tedrdg3}, \ref{08ujivxwfal879u}(ii) and (\ref{dimension-formula02}), there is $X_{21}\in V(G_{p^s}(n,1))$ such that
$\scriptsize X_2=\left(\begin{array}{c}
 X_{21} \\
 X_{22} \\
 \end{array}
 \right)$ and $\scriptsize X=\left(\begin{array}{c}
 X_1 \\
 X_{22} \\
 \end{array}
 \right)$, where $X_{22}=X_1\cap X_2$ is a $2$-subspace. Note that $\varphi_{m-2}(X_1)\stackrel{mc}{\sim} \varphi_{m-2}(X_2)$.
 Applying (\ref{automorphism021}) and (\ref{automorphism026}), we have
$$\varphi_{m-3}(X)=\varphi_{m-2}(X_1)\vee\varphi_{m-2}(X_2)=X_1U\vee X_2U=\left(\begin{array}{c}
                                                X_1 \\
                                                 X_{22} \\
                                              \end{array}
                                            \right)U =XU.$$
 Therefore, we get
\begin{equation}\label{automorphism027}
\mbox{$\varphi_{m-3}(X)=XU$,  \ $X\in V(G_{p^s}(n,3))$.}
\end{equation}

Similarly,  we can prove that
\begin{equation}\label{automorphism028}
\mbox{$\varphi(X)=\varphi_0(X)=XU$,  \ $X\in V(G_{p^s}(n, m))$.}
\end{equation}

{\bf Case 2}. \ $\varphi$ maps stars to tops with $n=2m$.

Let $\psi(X)=(\varphi(X))^\perp$,  $X\in V(G_{p^s}(2m,m))$.
By (\ref{Dualsubspace02}) and (\ref{Dualsubspace04}),  $\psi$ is surjection from $V(G_{p^s}(2m,m))$ to itself.  By (\ref{arithmeticdistance07}),
$\psi$ is also an automorphism of $G_{p^s}(2m,m)$. Using (\ref{Grassmann009bb}), $\psi$ maps  stars to  stars. By the Case 1,
there exists  $U\in GL_n(\mathbb{Z}_{p^s})$ such that $\psi(X)=XU$, $X\in V(G_{p^s}(2m, m))$. Thus
\begin{equation}\label{automorphism030}
\mbox{$\varphi(X)=(XU)^\perp $,  \ $X\in V(G_{p^s}(2m, m))$.}
\end{equation}

By above discussion, we complete the proof of Theorem \ref{automorphism014}.
 \hfill          $\Box$

\begin{rem}\label{automorphism029} \
Suppose that $2\leq m<n< 2m$. Then every  maximum clique of $G_{p^s}(n,m)$ is a top. Recall that $G_{p^s}(n,m)$ is isomorphic to $G_{p^s}(n,n-m)$.
By the bijection $Y\mapsto Y^\perp$, it is easy to see that every automorphism $\varphi$ of $G_{p^s}(n,m)$
is of the form $\varphi(X)=(XU)^\perp$, where $U\in GL_n(\mathbb{Z}_{p^s})$.
\end{rem}



\end{document}